\documentclass[a4paper,12pt]{amsart}

\usepackage{amsmath,amssymb}
\usepackage{graphicx}
\usepackage[dvips,arrow]{xy}
\usepackage{pb-diagram,pb-xy}

\textwidth=\paperwidth \advance\textwidth by-2\oddsidemargin
\advance\oddsidemargin by-1in
\evensidemargin=\oddsidemargin
\calclayout

\def\F{\mathcal{F}}
\def\G{\mathcal{G}}
\def\H{\mathcal{H}}
\def\K{\mathcal{K}}

\def\Z{\mathbb{Z}}
\def\Q{\mathbb{Q}}

\def\hat{\widehat}
\def\tilde{\widetilde}

\def\id{\mathrm{id}}

\def\^#1{^{(#1)}}

\theoremstyle{plain}
\newtheorem{theorem}{Theorem}[section]
\newtheorem{proposition}[theorem]{Proposition}
\newtheorem{corollary}[theorem]{Corollary}
\newtheorem{lemma}[theorem]{Lemma}

\theoremstyle{definition}
\newtheorem{definition}[theorem]{Definition}

\newtheorem{remark}[theorem]{Remark}

\def\to{\mathchoice{\longrightarrow}{\rightarrow}{\rightarrow}{\rightarrow}}
\makeatletter
\newcommand{\shortxra}[2][]{\ext@arrow 0359\rightarrowfill@{#1}{#2}}
\def\longrightarrowfill@{\arrowfill@\relbar\relbar\longrightarrow}
\newcommand{\longxra}[2][]{\ext@arrow 0359\longrightarrowfill@{#1}{#2}}
\renewcommand{\xrightarrow}[2][]{\mathchoice{\longxra[#1]{#2}}%
  {\shortxra[#1]{#2}}{\shortxra[#1]{#2}}{\shortxra[#1]{#2}}}
\makeatother

\begin{document}

\title [Injectivity theorems and algebraic closures]%
{Injectivity theorems and algebraic closures of groups with
  coefficients}

\author{Jae Choon Cha}

\address{Information and Communications University, Daejeon 305--732,
  Korea}

\email{jccha@icu.ac.kr}

\def\subjclassname{\textup{2000} Mathematics Subject Classification}
\expandafter\let\csname subjclassname@1991\endcsname=\subjclassname
\expandafter\let\csname subjclassname@2000\endcsname=\subjclassname


\begin{abstract}
  Recently, Cochran and Harvey defined torsion-free derived series of
  groups and proved an injectivity theorem on the associated
  torsion-free quotients.  We show that there is a universal
  construction which extends such an injectivity theorem to an
  isomorphism theorem.  Our result relates injectivity theorems to a
  certain homology localization of groups.  In order to give a
  concrete combinatorial description and existence proof of the
  necessary homology localization, we introduce a new version of
  algebraic closures of groups with coefficients by considering a
  certain type of equations.

\end{abstract}

\maketitle

\section{Introduction}
\label{section:introduction}

Let $R$ be a subring of the rationals.  A map $f\colon X \to Y$
between two spaces $X$ and $Y$ is called an \emph{$R$-homology
  equivalence} if $f$ induces isomorphisms on $H_*(-;R)$.  Homology
equivalences play an important role in the study of various problems
in geometric topology, including homology cobordism of manifolds and
concordance of embeddings, in particular knot and link concordance.
In this regard, a central question is how the fundamental groups of
homology equivalent spaces relate.  As a preliminary observation, it
can be easily seen that an $R$-homology equivalence induces a
homomorphism on the fundamental groups which is 2-connected on
$R$-homology; we say that a group homomorphism $\phi$ is
\emph{2-connected on $R$-homology} if $\phi$ induces an isomorphism on
$H_1(-;R)$ and a surjection on $H_2(-;R)$.  When $R=\Z$ (resp.\ $\Q$),
$\phi$ is called \emph{integrally (resp.\ rationally) 2-connected}.
We also remark that in most applications it suffices to consider
finite complexes (for example, compact manifolds) and so fundamental
groups can be assumed to be finitely presented.

Probably the first landmark result on the relationship between
homology equivalences and fundamental groups is an isomorphism theorem
of Stallings, which says that an integrally 2-connected homomorphism
$\pi\to G$ induces an isomorphism $\pi/\pi_q \to G/G_q$, where $G_q$
denotes the $q$-th lower central subgroup of $G$
\cite{Stallings:1965-1}.  This has several well-known topological
applications: abelian invariants such as linking numbers can be viewed
as an application of the simplest nontrivial case $G/G_2=H_1(G)$ of
Stallings' theorem.  More generally, the invariance of $G/G_q$ plays a
key role in understanding Milnor's $\bar\mu$-invariants of links
\cite{Milnor:1954-1,Milnor:1957-1} as shown in Casson's
work~\cite{Casson:1975-1}.  Orr defined further homotopy invariants of
links using Stallings' theorem~\cite{Orr:1987-1,Orr:1989-1}.  In some
cases homology cobordism invariants can be obtained by combining
Stallings' theorem and the index theorem in a similar way as Levine's
work on Atiyah--Patodi--Singer signatures of
links~\cite{Levine:1994-1}.  In \cite{Friedl:2003-1}, Friedl applied
this method to reformulate and generalize a link concordance invariant
obtained from certain nonabelian and irregular covers due to the
author and Ko~\cite{Cha-Ko:1999-1}.

In 2004, Cochran and Harvey announced a remarkable discovery of an
injectivity theorem relating the rational homology of a group $G$ to a
certain type of derived series $\{G_H\^n\}$, which is called the
\emph{torsion-free derived series}~\cite{Cochran-Harvey:2004-1}.  (The
series is due to Harvey~\cite{Harvey:2005-1}; for a precise definition
of $G_H\^n$, see \cite{Cochran-Harvey:2004-1,Harvey:2005-1}, or
Section~\ref{section:torsion-free-derived-series}.)  The main result
of~\cite{Cochran-Harvey:2004-1} can be stated as follows: if $\pi \to
G$ is a rationally 2-connected homomorphism of a finitely generated
group $\pi$ into a finitely presented group $G$, then it induces an
injection $\pi/\pi_H\^n \to G/G_H\^n$.  It has interesting
applications as illustrated in a recent result of Harvey; she obtained
invariants of homology cobordism by combining the injectivity with
rank invariants and $L^{(2)}$-signature invariants.

The advent of the injectivity theorem leads us to ask a natural
question: can one extend the torsion-free quotient $G/G_H\^n$ in such
a way that an isomorphism is induced instead of an injection?  More
generally, when can such an injectivity theorem be extended to an
isomorphism theorem?  Regarding topological applications, we remark
that the fundamental ideas of the applications of Stallings' theorem
could be reused when one has an isomorphism theorem.

In this paper we study injectivity theorems and their extensions to
isomorphism theorems in a general setting motivated from
Cochran--Harvey's result.  To formalize injectivity theorems, we
introduce a notion of an \emph{I-functor}.  For an arbitrary
coefficient ring $R\subset \Q$, we think of the collection $\Omega^R$
of homomorphisms of finitely generated groups into finitely presented
groups which are 2-connected on $R$-homology.  Roughly speaking, we
define an \emph{I-functor} $\H$ to be a functorial association $G \to
\H(G)$ such that to each $\pi\to G$ in $\Omega^R$, an injection
$\H(\pi)\to \H(G)$ is associated.  Of course the main example of an
I-functor is $\H(G)=G/G_H\^n$ where $R=\Q$.  We also formalize an
isomorphism theorem extending the injectivity as follows: a
\emph{container} of an I-functor $\H$ is defined to be another
I-functor $\F$ such that $\H$ injects into $\F$, that is,
$\H(G)\subset \F(G)$, and $\F$ associates an isomorphism to each
morphism in $\Omega^R$.  (Because there is some technical
sophistication, we postpone precise definitions to
Section~\ref{section:I-functors-and-containers}; here we just remark
that everything is required to have certain functorial properties
which are naturally expected.)

Then our question can be stated as whether there is a container of a
given I-functor.  Note that if one has a container, then it is easy to
construct a larger container by extending it; for example, take the
direct sum with a constant functor.  So the most essential one is a
minimal container.  We prove the following result:

\begin{theorem}
  \label{theorem:rough-version-of-main-theorem}
  If an I-functor $\H$ commutes with limits, then there exists a
  container $\hat\H$ of $\H$ which is universal (initial) in the
  following sense: if $\F$ is another container, then $\hat\H$ injects
  into $\F$ in a unique way.
\end{theorem}

In other words, $\hat\H$ provides an isomorphism theorem which extends
the injectivity of $\H$, and it is universal among such extensions.
For a more precise statement, see
Section~\ref{section:I-functors-and-containers}.  We remark that from
its universality it follows that $\hat\H$ is a unique minimal
container of~$\H$.

As a corollary of Theorem~\ref{theorem:rough-version-of-main-theorem},
we show that the torsion-free derived quotient $G \to G/G_H\^n$ has a
universal container (see
Corollary~\ref{corollary:existence-of-container-of-C-H-quotient}).  We
remark that this special case was partially addressed
in~\cite{Cochran-Harvey:2004-1}; they constructed a container of
$G/G_H\^n$ by using an iterated semidirect product of certain homology
groups.  We also show that the container in
\cite{Cochran-Harvey:2004-1} fails to be universal in our sense (see
Theorem~\ref{theorem:non-universality-of-C-H-container}).  An
interesting observation is that this is closely related to the use of
the \emph{Ore localization} of a group ring $\Q[G]$ of a
poly-torsion-free-abelian group $G$ in the construction of the
container in~\cite{Cochran-Harvey:2004-1}.  In the Ore localization
all nonzero elements are inverted, but it turns out that this is too
excessive; in showing that our universal container is strictly
smaller, it is illustrated that the unnecessarily inverted elements
are ones in the kernel of the augmentation $\Q[G] \to \Q$.  (See
Section~\ref{section:torsion-free-derived-series} for more details.)
This gives us a motivation for expecting a similar but more natural
theory using the \emph{Cohn localization}, instead of the Ore
localization.

On the other hand, in our results there is something beyond the
existence of a universal container.  It relates injectivity theorems
to a certain localization functor of groups.  In general, one can view
a localization functor as a \emph{universal} construction inverting a
given collection $\Omega$ of morphisms in a category.  (Our definition
of a localization is given in
Section~\ref{section:localization-wrt-2-connected-morphisms}.)  So, if
a localization with respect to $\Omega$ exists, it provides a natural
isomorphism theorem for morphisms in $\Omega$.  That is, it associates
an isomorphism to each morphism in~$\Omega$.  In order to prove
Theorem~\ref{theorem:rough-version-of-main-theorem}, we use a
particular localization functor $G \to \hat G$ with respect to the
collection $\Omega^R$ considered above.  In fact, in the proof of
Theorem~\ref{theorem:rough-version-of-main-theorem}, we show that the
universal container $\hat\H$ of $\H$ is given by $\hat\H(G)=\H(\hat
G)$.  This presents another point of view that the injective
homomorphism induced by $\H$ can be regarded as a restriction of a
natural isomorphism obtained from the localization; for any
homomorphism in $\Omega^R$, the induced isomorphism on the
localization gives rise to an isomorphism on $\hat\H(-)$, and the
induced homomorphism on $\H(-)$ is its restriction.

We remark that several homology localizations of groups have been
studied in the literature, including works of Bousfield
\cite{Bousfield:1975-1}, Vogel, Levine
\cite{Levine:1989-2,Levine:1989-1}, and
Farjoun--Orr--Shelah~\cite{Farjoun-Orr-Shelah:1989}.  In particular,
in \cite{Levine:1989-2,Levine:1989-1} Levine introduced the notion of
an algebraic closure of a group, which turns out to be equivalent to
Vogel's localization with respect to integrally 2-connected
homomorphisms from finitely generated groups into finitely presented
groups which are normally surjective.  The localization with respect
to our $\Omega^R$ which is used to prove
Theorem~\ref{theorem:rough-version-of-main-theorem} is similar to that
of Levine, but distinguished in two points: first the normal
surjectivity condition is not required, and second, an arbitrary
subring $R$ of $\Q$ is used as homology coefficients.  Although it
may be regarded as folklore that there exists a localization with
respect to $\Omega^R$, we give an existence proof for concreteness
since we could not find any published one in the literature.  In
addition our work provides a combinatorial description of the desired
$R$-homology localization.  For this purpose, we introduce a new
version of algebraic closures of groups, modifying the idea of
Farjoun, Orr, and Shelah~\cite{Farjoun-Orr-Shelah:1989} and
Levine~\cite{Levine:1989-2,Levine:1989-1}.  We think of a certain
class of systems of equations over a group $G$ which we call
\emph{$R$-nullhomologous}, and define an \emph{(algebraic)
  $R$-closure} $\hat G$ of $G$ in terms of solubility of such systems
of equations.  We show the existence of an $R$-closure $\hat G$ for
any group $G$, and show that it is equivalent to the desired
localization with respect to~$\Omega^R$.

The remaining part of this paper is organized as follows.  In
Section~\ref{section:I-functors-and-containers} we prove the existence
of a universal container, assuming the existence of the $R$-closure of
a group.  In Section~\ref{section:torsion-free-derived-series} we
apply the results of Section~\ref{section:I-functors-and-containers}
to the case of the torsion-free derived quotient, and show that our
universal container is strictly smaller than a container constructed
in~\cite{Cochran-Harvey:2004-1}.  In the remaining sections, we prove
results on the $R$-closure which are used in previous sections.  In
Section~\ref{section:nullhomologous-equations} we introduce the notion
of $R$-nullhomologous systems of equations, and in
Section~\ref{section:localization-wrt-2-connected-morphisms} it is
related to the localization of groups with respect to~$\Omega^R$.  In
Section~\ref{section:existence-of-closures} we show the existence of
the $R$-closure and some properties of the $R$-closure of a finitely
presented group.

\subsection*{Acknowledgements}

The author would like to thank Stefan Friedl and Kent Orr for
discussions from which the fundamental idea of this paper is obtained.
The author also thanks Tim Cochran and Shelly Harvey who kindly
provided a copy of slides containing their results in
\cite{Cochran-Harvey:2004-1,Harvey:2005-1} before the manuscripts
became available.  Finally, comments from an anonymous referee were
very helpful in improving this paper.

\section{I-functors and containers}
\label{section:I-functors-and-containers}

We start with a formalization of the notion of a container of the
torsion-free derived quotient $G/G_H\^n$.  Here we have a technical
issue that the association $G \to G/G_H\^n$ is not a functor of the
category $\G$ of groups; not all group homomorphisms induce a morphism
on the quotients, although the result of Cochran--Harvey guarantees
that a rationally 2-connected homomorphism of a finitely generated
group into a finitely presented group gives rise to an induced
homomorphism.

This leads us to consider what follows.  Recall that for any subring
$R$ of $\Q$, we denote by $\Omega^R$ the class of group homomorphisms
$\alpha\colon \pi \to G$ such that $\pi$ is finitely generated, $G$ is
finitely presented, and $\alpha$ is 2-connected on $R$-homology, that
is, $\alpha$ induces an isomorphism on $H_1(-;R)$ and a surjection on
$H_2(-;R)$.  Denoting $\H(G) = G/G_H\^n$, in
\cite{Cochran-Harvey:2004-1} it was shown that $\H$ has the following
properties for $R=\Q$:

\begin{enumerate}
\item To each group $G$, a homomorphism $p_G \colon G \to \H(G)$ is
  associated.
\item To each homomorphism $\alpha \colon \pi \to G$ in $\Omega^R$, an
  injection $\H(\pi) \to \H(G)$ is associated.
\item The above associations have naturality, that is,
  $\H(\beta\circ\alpha) = \H(\beta)\circ \H(\alpha)$, $\H(\alpha)\circ
  p_\pi = p_G \circ \alpha$ (see the diagrams below), and
  $\H(\id_G)=\id_{\H(G)}$ whenever the involved homomorphisms exist.
  \[
  \begin{diagram}
    \node{\H(\pi)} \arrow[2]{e,t}{\H(\beta\circ\alpha)} \arrow{se,b}{\H(\alpha)}
    \node[2]{\H(P)}
    \\
    \node[2]{\H(G)} \arrow{ne,b}{\H(\beta)}
  \end{diagram}
  \qquad
  \begin{diagram}
    \node{\pi}\arrow{e,t}{\alpha} \arrow{s,l}{p_\pi}
    \node{G} \arrow{s,r}{p_G}
    \\
    \node{\H(\pi)} \arrow{e,b}{\H(\alpha)}
    \node{\H(G)}
  \end{diagram}
  \]
  \xdef\lastcount{\arabic{enumi}}
\end{enumerate}

While $\H$ is not a functor of the category $\G$ of groups, we can
view $\H$ as a functor of a subcategory of finitely presented groups,
which are of our main interest regarding topological applications: let
$\G^R$ be the category whose objects are finitely presented groups and
whose morphisms are homomorphisms between finitely presented groups
which are 2-connected on $R$-homology.  Then from the above properties
it follows that $\H$ induces a functor $\G^R \to \G$.  Also $p$
induces a natural transformation from the obvious inclusion functor
$\G^R \to \G$ to (the functor induced by)~$\H$.  In particular,
homology equivalences between finite complexes give rise to morphisms
in $\G^R$ and then one can apply $\H$ to obtain injections.

One more obvious property of $\H$ which might be easily overlooked is
the following:
\begin{enumerate}
  \setcounter{enumi}{\lastcount}
\item To any isomorphism $\alpha\colon \pi\to G$, an isomorphism
  $\H(\alpha)\colon \H(\pi) \to \H(G)$ is associated.
\end{enumerate}
We remark that (4) does not follow from (1)--(3) since (1)--(3) do not
guarantee that $\H(\alpha)$ is defined for an isomorphism~$\alpha$ in
general.

Results of this section are not specific to the torsion-free
quotients; we consider any association $\H$ with the properties
above.

\begin{definition}
  $(\H,p)$ is called an \emph{I-functor with respect to
    $R$-coefficients} if the above (1)--(4) are satisfied.
\end{definition}

Note that if $R'$ is a subring of $R$, then an $I$-functor with
respect to $R$-coefficients is automatically an $I$-functor with
respect to $R'$-coefficients.  When the coefficient ring $R$ is
clearly understood, we simply say that $\H$ is an I-functor.

I-functors form (objects of) a category; a morphism $\tau$ between two
I-functors $\H$ and $\H'$ is defined to be a natural transformation
$\tau \colon \H \to \H'$, where $\H$ and $\H'$ are viewed as functors
$\G^R \to \G$ as an abuse of notation, such that the diagram
\[
\begin{diagram}
  \node{G}\arrow{s,l}{p_G}\arrow{se,t}{p'_G}\\
  \node{\H(G)}\arrow{e,b}{\tau_G} \node{\H'(G)}
\end{diagram}
\]
commute for each object $G$ in $\G^R$ where $p_G$ and $p'_G$ are the
natural transformations that the I-functors $\H$ and $\H'$ are endowed
with, respectively.  If each $\tau_G$ is injective, then we say that
$\tau$ is \emph{injective} and $\H'$ is an \emph{extension} of~$\H$.

\begin{definition}
  For an I-functor $\H$, a morphism $\tau\colon \H \to \F$ into
  another I-functor $\F$ is called a \emph{container} of $\H$ if
  $\tau$ is injective and $\F(\alpha)\colon \F(\pi)\to \F(G)$ is an
  isomorphism for any morphism $\alpha\colon \pi\to G$ in~$\Omega^R$.
\end{definition}

Sometimes we say that $\F$ is a container of $\H$ when we do not have
to specify $\tau\colon \H \to \F$ explicitly.

As mentioned in the introduction, we are interested in a universal
(initial) container of a given I-functor~$\H$.  To give its
definition, we consider the category of containers and injective
morphisms; objects are containers $\F$ of $\H$, and morphisms from
$\F$ to $\F'$ are injective morphisms $\F \to \F'$ between the two
I-functors $\F$ and $\F'$ which makes the diagram
\[
\begin{diagram}
  \node{\H}\arrow{s}\arrow{se}\\
  \node{\F}\arrow{e}
  \node{\F'}
\end{diagram}
\]
commute.

\begin{definition}
  A universal (initial) object $\F$ in the category of containers of
  $\H$ is called a \emph{universal} container of~$\H$, that is, for
  any container $\F'$ of $\H$, there is a unique morphism from $\F$
  to~$\F'$.
\end{definition}

Obviously a universal container is unique if it exists.  Also, a
universal container is automatically minimal, in the sense that it is
not a proper extension of another container.  So if a universal
container exists, it is a unique minimal container.

For our existence result of a universal container, we need to
formulate a relationship of an I-functor $\H$ and limits.  In this
paper it suffices to consider the direct limit of a sequence
\[
G_0 \to G_1 \to G_2 \to \cdots.
\]
of group homomorphisms in $\G^R$.
Usually, if $\H$ were an ordinary functor $\G\to \G$, we would say
that $\H$ commutes with limits when $\varinjlim \H(G_k) \cong
\H(\varinjlim G_k)$; more precisely, the isomorphism is explicitly
specified in this case.  Namely, $G_k \to \varinjlim G_k$ induces
$\H(G_k) \to \H(\varinjlim G_k)$, and then
\[
\varinjlim \H(G_k) \to \H(\varinjlim G_k)
\]
is induced.  If it is an isomorphism, then we say that $\H$ commutes
with limits.  However, in our case, because $\H$ is just an I-functor,
the homomorphism $G_k \to \varinjlim G_k$ does not necessarily induce
$\H(G_k) \to \H(\varinjlim G_k)$ in general.  So we need to adopt the
existence of this induced homomorphism as a part of a definition:

\begin{definition}
  An I-functor $\H$ is said to \emph{commute with limits} if for any
  sequence
  \[
  G_0 \to G_1 \to G_2 \to \cdots
  \]
  of morphisms $G_k \to G_{k+1}$ in $\G^R$, $\H$ associates to $G_k
  \to \varinjlim G_k$ a homomorphism $\H(G_k) \to \H(\varinjlim G_k)$,
  and its limit
  \[
  \varinjlim \H(G_k) \to \H(\varinjlim G_k)
  \]
  is an isomorphism.
\end{definition}

We note that even though each $G_k \to G_{k+1}$ is in $\G^R$,
$\varinjlim G_k$ is not necessarily (an object) in~$\G^R$.

Now we can precisely state the main result of this section.

\begin{theorem}
  \label{theorem:existence-of-container-of-I-functor}
  Suppose that $R$ is a subring of $\Q$ and $\H$ is an I-functor with
  respect to $R$-coefficients which commutes with limits.  Then there
  exists a universal container $\tau \colon \H \to \hat\H$ of $\H$,
  that is, for any container $\sigma\colon \H \to \F$, there is a
  unique injective morphism $\hat\sigma \colon \hat\H \to \F$ such
  that the diagram
  \[
  \begin{diagram}
    \node[2]{G}\arrow{ssw}\arrow{s}\arrow{sse}
    \\
    \node[2]{\hat \H(G)}\arrow{se,b,..}{\hat\sigma_G}
    \\
    \node{\H(G)} \arrow{ne,b}{\tau_G} \arrow[2]{e,b}{\sigma_G}
    \node[2]{\F(G)}
  \end{diagram}
  \]
  commutes.
\end{theorem}

To prove Theorem~\ref{theorem:existence-of-container-of-I-functor}, we
use a homology localization functor $E\colon \G \to \G$ with respect
to $R$-coefficients.  At this moment we just need the following
properties of~$E$, which are analogues of Levine's results on
algebraic closures of groups \cite{Levine:1989-2,Levine:1989-1}; a
construction of our $E$ and proofs of the necessary properties are
postponed to later sections.

\begin{theorem}
  \label{theorem:properties-of-closures}
  For any subring $R$ of $\Q$, there is a pair $(E,i)$ of a functor
  $E\colon \G \to \G$ and a natural transformation $i \colon \id_\G
  \to E$ which has the following properties:
  \begin{enumerate}
  \item For any $\alpha\colon \pi \to G$ in $\Omega^R$, the induced
    homomorphism $E(\alpha) \colon E(\pi) \to E(G)$ is an isomorphism.
  \item For any object $G$ in $\G^R$, there is a sequence
    \[
    G=G_0 \to G_1 \to \cdots \to G_k \to \cdots
    \]
    of morphisms $G_k \to G_{k+1}$ in $\G^R$ such that $E(G)=
    \varinjlim G_k$ and $i_G \colon G \to E(G)$ is the limit
    homomorphism.
  \end{enumerate}
\end{theorem}

We denote $E(G)$ by $\hat G$.  For any I-functor $\H$ which commutes
with limits, we will prove that the composition of $E$ and $\H$ is a
universal container of $\H$.

\begin{proof}
  [Proof of Theorem~\ref{theorem:existence-of-container-of-I-functor}]
  We define $\hat\H(G) = \H(\hat G)$ and $\hat p_G \colon G \to
  \hat\H(G)$ to be the composition
  \[
  G \xrightarrow{i_G} \hat G \xrightarrow{p_{\hat G}} \H(\hat G)=\hat\H(G).
  \]
  In other words, $\hat\H = \H \circ E$ and $\hat p = p \circ i$.  We
  will show that $(\hat\H,\hat p)$ is an I-functor.  For any
  $\alpha\colon \pi \to G$ which is in $\Omega^R$, $\hat\alpha\colon
  \hat\pi \to \hat G$ is an isomorphism.  Applying $\H$, we obtain an
  induced isomorphism $\H(\hat\pi) \to \H(\hat G)$.  We define
  $\hat\H(\alpha)$ to be this isomorphism.  Viewing $(\hat H, \hat p)$
  as $(\H \circ E, p \circ i)$, the required naturality of $(\hat H,
  \hat p)$ follows from that of $(\H,p)$ and $(E,i)$.

  For a finitely presented group $G$, $(\hat\H, \hat p)$ can be
  interpreted as follows.  Choose a sequence
  \[
  G=G_0 \to G_1 \to G_2 \to \cdots
  \]
  of morphisms in $\G^R$ such that $i_G \colon G \to \hat G$ is the
  limit homomorphism $G \to \varinjlim G_k \cong \hat G$.  By the naturality
  of $p$,
  \[
  \begin{diagram}
    \node{G} \arrow{e} \arrow{s,l}{p_G}
    \node{G_k} \arrow{e} \arrow{s,r}{p_{G_k}}
    \node{\varinjlim G_k} \arrow{s,r}{p_{\varinjlim G_k}}
    \\
    \node{\H(G)} \arrow{e}
    \node{\H(G_k)} \arrow{e}
    \node{\H(\varinjlim G_k)}
  \end{diagram}
  \]
  commutes.  Taking the limit, we have a commutative diagram
  \[
  \begin{diagram}
    \node{G} \arrow{e,t}{i_G} \arrow{s,l}{p_G}
    \node{\varinjlim G_k} \arrow{e,=} \arrow{s,r}{\varinjlim p_{G_k}}
    \node{\varinjlim G_k \hbox to 0mm{ $= \hat G$\hss}}
      \arrow{s,r}{p_{\varinjlim G_k}}
    \\
    \node{\H(G)} \arrow{e,b}{\text{limit map}}
    \node{\varinjlim \H(G_k)} \arrow{e,b}{\cong}
    \node{\H(\varinjlim G_k)\hbox to 0mm{ $= \H(\hat G)$\hss}}
  \end{diagram}
  \hphantom{= \H(\hat G)}
  \]
  That is, $\hat\H(G) = \varinjlim \H(G_k)$ and $\H$ associates to
  $i_G$ the limit homomorphism
  \[
  \H(i_G)\colon \H(G) \to \hat\H(G) = \varinjlim \H(G_k).
  \]
  Also, $\hat p_G\colon G \to \hat\H(G)$ is the composition
  \[
  G \xrightarrow{p_G} \H(G) \xrightarrow{\H(i_G)} \varinjlim \H(G_k).
  \]

  Now we construct an injective morphism $\tau\colon (\H,p) \to
  (\hat\H, \hat p)$ between the I-functors $(\H,p)$ and $(\hat\H, \hat
  p)$ as follows.  For a finitely presented group $G$, there exists
  $\H(i_G)\colon \H(G) \to \H(\hat G)$ as discussed above.  We define
  $\tau_G \colon \H(G) \to \hat\H(G)$ to be $\H(i_G)$, that is, $\tau
  = \H \circ i$.  Viewing $\tau$ as a transformation between functors
  $\H, \H'\colon \G^R \to \G$, the naturality of $\tau$ follows from
  that of $\H$ and~$i$.  Furthermore,
  \[
  \begin{diagram}
    \node{G} \arrow{s,l}{p_G} \arrow{se,t}{\hat p_G}
    \\
    \node{\H(G)} \arrow{e,b}{\tau_G}
    \node{\hat\H(G)}
  \end{diagram}
  \]
  commutes since $\hat p_G = \H(i_G) \circ p_G$.  This shows that
  $\tau$ is a morphism $(\H,p) \to (\hat\H,\hat p)$.  To show the
  injectivity, we consider a sequence $G=G_0 \to G_1 \to \cdots$ with
  limit $\hat G$ as above.  Since $\H$ is an I-functor and $G \to G_k$
  is in $\Omega^R$, $\H(G) \to \H(G_k)$ is injective.  Since $\tau_G =
  \H(i_G)$ is the limit of $\H(G) \to \H(G_k)$, $\tau_G$ is injective
  too.

  We will show that $\tau\colon (\H,p) \to (\hat\H,\hat p)$ has the
  universal property.  Suppose that $\sigma\colon (\H,p) \to (\F,q)$
  is a container.  We define a morphism $\hat\sigma \colon
  (\hat\H,\hat p) \to (\F,q)$ as follows: for a finitely presented
  group $G$, choose $G=G_0 \to G_1 \to \cdots$ with limit $\hat G$ as
  above.  Taking the limit of
  \[
  \begin{diagram}\dgARROWLENGTH=1.9em
      \node[2]{G} \arrow{e}\arrow{s,r}{p_G}\arrow{ssw,t}{q_G}
      \node{G_k} \arrow{s,l}{p_{G_k}} \arrow{sse,t}{q_{G_k}}
      \\
      \node[2]{\H(G)} \arrow{e,J} \arrow{sw,b,J}{\sigma_G}
      \node{\H(G_k)} \arrow{se,b,J}{\sigma_{G_k}}
      \\
      \node{\F(G)} \arrow[3]{e,b}{\cong}
      \node[3]{\F(G_k)}
    \end{diagram}
  \]
  we obtain a commutative diagram
  \begin{equation}
    \begin{diagram}\dgARROWLENGTH=1.9em
      \node[2]{G} \arrow{e}\arrow{s,r}{p_G}\arrow{ssw,t}{q_G} 
      \node{\varinjlim G_k}
      \arrow{s,l}{\varinjlim p_{G_k}} \arrow{sse,t}{\varinjlim q_{G_k}}
      \\
      \node[2]{\H(G)} \arrow{e,J} \arrow{sw,b,J}{\sigma_G}
      \node{\varinjlim \H(G_k)} \arrow{se,b,J}{\varinjlim \sigma_{G_k}}
      \\
      \node{\F(G)} \arrow[3]{e,b}{\cong}
      \node[3]{\varinjlim \F(G_k)}
    \end{diagram}
    \tag*{$(*)_{\{G_k\}}$}
  \end{equation}
  We define $\hat \sigma_G$ to be $\varinjlim \sigma_{G_k}$, that is,
  \[
  \hat\sigma_G \colon \hat\H(G) =\H(\hat G) = \varinjlim \H(G_k)
  \xrightarrow{\varinjlim \sigma_{G_k}} \varinjlim \F(G_k) = \F(G).
  \]
  Since each $\sigma_{G_k}$ is injective, so is~$\hat\sigma_G$.

  At the moment, our $\hat\sigma_G$ depends on the choice of
  $\{G_k\}$.  Before showing that it is well-defined, we prove the
  naturality of $\hat\sigma$.  Suppose $\alpha\colon \pi \to G$ is in
  $\G^R$ and $\pi=\pi_0 \to \pi_1 \to \cdots$ and $G=G_0 \to G_1 \to
  \cdots$ are sequences giving $\hat\pi$ and $\hat G$ as above.  First
  we consider a special case that $\{\pi_k\}$ and $\{G_k\}$ behave
  nicely under $\alpha$, that is, we suppose that there are
  homomorphisms $\pi_k \to G_k$ which fit into the following
  commutative diagram:
  \[
  \begin{diagram}
    \node{\hbox to 0mm{\hss$\pi=$ }\pi_0} \arrow{s,l}{\alpha}\arrow{e}
    \node{\pi_1} \arrow{s}\arrow{e}
    \node{\pi_2} \arrow{s}\arrow{e}
    \node{\cdots}
    \\
    \node{\hbox to 0mm{\hss$G=$ }G_0} \arrow{e}
    \node{G_1} \arrow{e}
    \node{G_2} \arrow{e}
    \node{\cdots}
  \end{diagram}
  \]
  Then $\alpha$, together with $\pi_k \to G_k$, induces a ``morphism''
  between the diagrams $(*)_{\{\pi_k\}}$ and $(*)_{\{G_k\}}$.  In
  particular, the morphisms $\hat\H(\alpha)$ and $\F(\alpha)$ give us
  the commutative diagram below, which says that $\hat\sigma$ is
  natural in this special case:
  \[
  \begin{diagram}
    \node{\hat\H(\pi)}\arrow{e,t}{\hat\H(\alpha)}
      \arrow{s,l}{\varinjlim\sigma_{\pi_k}=\hat\sigma_\pi} 
    \node{\hat\H(G)} \arrow{s,r}{\hat\sigma_G=\varinjlim\sigma_{G_k}}
    \\
    \node{\F(\pi)} \arrow{e,b}{\F(\alpha)}
    \node{\F(G)}
  \end{diagram}
  \]

  To reduce the general case to the special case above, we appeal to
  the following result, which is a horseshoe-type lemma:

  \begin{lemma}
    \label{lemma:horseshoe-special-case}
    Suppose that $\pi=\pi_0 \to \pi_1 \to \cdots$ and $G=G_0 \to G_1
    \to \cdots$ are sequences of morphisms in $\G^R$ such that
    $\hat\pi=\varinjlim \pi_k$ and $\hat G = \varinjlim G_k$, and
    $\alpha\colon \pi \to G$ is in $\G^R$.  Then there exists a
    commutative diagram
    \[
    \hphantom{G=}
    \begin{diagram}\dgARROWLENGTH=1em
      \node{\hbox to 0mm{\hss$\pi=$ }\pi_0} \arrow{e} \arrow{s,l}{\alpha}
      \node{\pi_1} \arrow{e} \arrow{s}
      \node{\cdots}  \arrow{e}
      \node{\pi_k} \arrow{e} \arrow{s}
      \node{\cdots}  \arrow{e}
      \node{\varinjlim \pi_k \hbox to 0mm{ $= \hat\pi$\hss}} \arrow{s}\\
      \node{\hbox to 0mm{\hss$G=$ }P_0} \arrow{e}
      \node{P_1} \arrow{e}
      \node{\cdots}  \arrow{e}
      \node{P_k} \arrow{e}
      \node{\cdots}  \arrow{e}
      \node{\varinjlim P_k} \\
      \node{\hbox to 0mm{\hss$G=$ }G_0} \arrow{e} \arrow{n,=}
      \node{G_1} \arrow{e} \arrow{n}
      \node{\cdots}  \arrow{e}
      \node{G_k} \arrow{e} \arrow{n}
      \node{\cdots}  \arrow{e}
      \node{\varinjlim G_k \hbox to 0mm{ $= \hat G$\hss}} \arrow{n}
    \end{diagram}\hphantom{= \hat G}
    \]
    where $\pi_k \to P_k$, $G_k \to P_k$, and $P_k \to P_{k+1}$ are in
    $\G^R$,
    \[
    \hat G =\varinjlim G_k \to \varinjlim P_k
    \]
    is an isomorphism, and the limit homomorphism
    \[
    G=P_0 \to \varinjlim P_k \cong \hat G
    \]
    is equal to $i_G \colon G \to \hat G$.
  \end{lemma}

  Applying the above special case to $(\{\pi_k\}, \{P_k\})$ and
  $(\{P_k\},\{G_k\})$, we obtain a commutative diagram
  \[
  \begin{diagram}
    \node{\hat\H(\pi)}\arrow{e,t}{\hat\H(\alpha)}\arrow{s,l}{\varinjlim\sigma_{\pi_k}} 
    \node{\hat\H(G)} \arrow{e,t}{\hat\H(\id)=\id}\arrow{s,r}{\varinjlim\sigma_{P_k}}
    \node{\hat\H(G)} \arrow{s,r}{\varinjlim\sigma_{G_k}}
    \\
    \node{\F(\pi)} \arrow{e,b}{\F(\alpha)}
    \node{\F(G)} \arrow{e,b}{\F(\id)=\id}
    \node{\F(G)}
  \end{diagram}
  \tag{$**$}
  \]
  This shows that $\hat\sigma$ behaves naturally for $\pi \to G$ even
  when $\hat\sigma_\pi$ and $\hat\sigma_G$ are defined using
  arbitrarily chosen $\{\pi_k\}$ and $\{G_k\}$.

  We can also use the same argument to show the well-definedness of
  $\hat\sigma_G$, that is, $\hat\sigma_G$ is independent of the choice
  of $\{G_k\}$.  For this, we apply
  Lemma~\ref{lemma:horseshoe-special-case} to a special case that
  $\pi=G$ and $\alpha\colon \pi \to G$ is the identity.  Then for any
  $\{\pi_k\}$ and $\{G_k\}$ with limit $\hat G$, there is $\{P_k\}$
  which gives the diagram~($**$).  Since $\hat\H(\alpha)=\id$ in this
  case, it follows that the homomorphisms
  \[
  \varinjlim \sigma_{\pi_k},\, \varinjlim
  \sigma_{P_k},\, \varinjlim \sigma_{G_k}\colon \hat\H(G) \to \F(G)
  \]
  are all equal.
  
  From the diagram~$(*)_{\{G_k\}}$, we obtain a commutative diagram
  \[
  \begin{diagram}
    \node[2]{G} \arrow{ssw,t}{\hat p_G} \arrow{sse,t}{q_G} \arrow{s,l,3}{p_G}
    \\
    \node[2]{\H(G)} \arrow{sw,b}{\tau_G} \arrow{se,b}{\sigma_G}
    \\
    \node{\hat\H(G)} \arrow[2]{e,b}{\hat\sigma_G}
    \node[2]{\F(G)}
  \end{diagram}.
  \]
  From this it follows that $\hat\sigma$ can be viewed as a morphism
  between containers.

  Finally we show the uniqueness of $\hat\sigma$.  Suppose that
  $\hat\sigma' \colon (\hat\H,\hat p) \to (\F,q)$ is another morphism
  between the containers $(\hat\H,\hat p)$ and $(\F,q)$.  For a
  sequence $G=G_0 \to G_1 \to \cdots$ giving $\hat G$, we have the
  following commutative diagram:
  \[
  \begin{diagram}\dgARROWLENGTH=.7em
    \node[3]{\H(G)}
      \arrow[2]{sw,t}{\tau_G} \arrow[2]{s,r}{\sigma_G}
      \arrow[2]{se}
    \\
    \\
    \node{\hat\H(G)} \arrow[2]{e,t}{\hat\sigma'_G} \arrow[2]{se,b}{\cong}
    \node[2]{\F(G)} \arrow[2]{se,t,1}{\cong}
    \node[2]{\H(G_k)} \arrow{sw,-} \arrow[2]{s,r}{\sigma_{G_k}}
    \\
    \node[4]{}\arrow{sw,t}{\tau_{G_k}}
    \\
    \node[3]{\hbox to0mm{\hss$\varinjlim\H(G_k)=$ }\hat\H(G_k)}
      \arrow[2]{e,b}{\hat\sigma'_{G_k}}
    \node[2]{\F(G_k)}
  \end{diagram}
  \]
  Since $\hat\sigma'_{G_k}\tau_{G_k} = \sigma_{G_k} = \tau_{G_k}
  \hat\sigma_{G_k} $, $\hat\sigma'_{G}$ and $\hat\sigma_{G}\colon
  \hat\H(G) \to \F(G)$ coincide on the image of $\tau_{G_k}\colon
  \H(G_k) \to \hat\H(G_k)=\hat\H(G)$.  From this it follows that
  $\hat\sigma'_{G}=\hat\sigma_{G}$ since
  $\hat\H(G)=\varinjlim\H(G_k)$.
\end{proof}

\section{Containers of torsion-free derived series}
\label{section:torsion-free-derived-series}

In this section we focus on a special case of an I-functor, namely the
torsion-free derived quotient $G \to G/G_H\^n$.  We begin by recalling
the definition of $G_H\^n$ in~\cite{Cochran-Harvey:2004-1}.  For a
group $G$, $G_H\^0$ is defined to be $G$ itself.  Suppose $G_H\^n$ has
been defined to be a normal subgroup of $G$ such that that $G/G_H\^n$
is a poly-torsion-free-abelian (PTFA) group.  Since the integral group
ring of a PTFA group is an Ore domain, there exists the skew field
$\K[G/G_H\^n]$ of (right) quotients of $\Z[G/G_H\^n]$, that is,
$\K[G/G_H\^n] = \Z[G/G_H\^n] (\Z[G/G_H\^n]-\{0\})^{-1}$.  Note that
$\K[G/G_H\^n]$ is $\Z[G/G_H\^n]$-flat.  $G_H\^{n+1}$ is defined to be
the kernel of the following composition:
\[
\begin{aligned}
  G_H\^n & \to G_H\^n/[G_H\^n,G_H\^n] = H_1(G;\Z[G/G_H\^n]) \\
  & \to H_1(G;\Z[G/G_H\^n])\mathbin{\mathop{\otimes}_{\Z[G/G_H\^n]}}
  \K[G/G_H\^n] = H_1(G;\K[G/G_H\^n])
\end{aligned}
\]
Since $G/G_H\^{n+1}$ is an extension of $G_H\^n/G_H\^{n+1}$ by
$G/G_H\^n$ and $G_H\^n/G_H\^{n+1}$ is a subgroup of
$H_1(G;\K[G/G_H\^n])$ which is a torsion-free abelian group,
$G/G_H\^{n+1}$ is PTFA so that one can continue this process.  For
$n=\omega$, the first infinite ordinal, $G_H\^n$ is defined to be the
intersection of $G_H\^k$ where $k$ runs over all integers.  For
further details see~\cite{Cochran-Harvey:2004-1}.

We denote $\H_n(G) = G/G_H\^n$.  It was shown in
\cite{Cochran-Harvey:2004-1} that $\H_n$, equipped with the projection
$G \to \H_n(G)$, is an I-functor with respect to $\Q$-coefficients.
(In this section the coefficient ring $R$ is always $\Q$.)  From
Theorem~\ref{theorem:existence-of-container-of-I-functor}, it follows
that $\H_n$ has a universal container:

\begin{corollary}
  \label{corollary:existence-of-container-of-C-H-quotient}
  There exists a universal container of~$\H_n$ for any~$n \le \omega$.
\end{corollary}

\begin{proof}
  By Theorem~\ref{theorem:existence-of-container-of-I-functor}, it
  suffices to show that $\H_n$ commutes with limits.  Suppose that
  \[
  G=G_0\to G_1 \to G_2 \to \cdots
  \]
  is a sequence of morphisms in $\G^R$.  We use the following two
  properties of~$\H_n$: first, $G \to \H_n(G)$ is obviously surjective
  for any $G$, and second, we need Lemma 5.3
  in~\cite{Cochran-Harvey:2004-1}: the natural map $G_k \to \varinjlim
  G_k$ gives rise to an injection $\H_n(G_k) \to \H_n(\varinjlim G_k)$
  for any $k$ and any $n\le \omega$.
    
  Taking the limit of the injections $\H_n(G_k) \to
  \H_n(\varinjlim G_k)$, we obtain an
  injection
  \[
  \varinjlim \H_n(G_k) \longrightarrow \H_n(\varinjlim G_k).
  \]
  Since every $x$ in $\H_n(\varinjlim G_k)$ is represented by an
  element in $\varinjlim G_k$, $x$ is in the image of some $G_k$, and
  so in the image of some $\H_n(G_k)$.  From this the surjectivity
  follows.
\end{proof}

Recall that the universal container $\hat\H_n$ of $\H_n$ is given by
$\hat\H_n(G)=\H_n(\hat G)$, where $\hat G$ denotes our homology
localization with respect to $\Q$-coefficients given in
Theorem~\ref{theorem:properties-of-closures}.

\begin{remark}
  For a $(4k-1)$-manifold $M$ with fundamental group $\pi$, Harvey
  considered the $L^{(2)}$-signature of $M$ associated to $\pi \to
  \H_n(\pi)$ as a homology cobordism invariant.  Since $\H_n(\pi) \to
  \hat \H_n(\pi) = \hat\pi/\hat\pi_H\^n$ is injective by
  Corollary~\ref{corollary:existence-of-container-of-C-H-quotient},
  from the induction property of the $L^{(2)}$-signature it follows
  that Harvey's invariant coincides with the $L^{(2)}$-signature
  associated to a characteristic quotient of the localization of
  $\pi$, namely $\pi \to \hat \pi \to \hat\pi/\hat\pi_H\^n$.
\end{remark}

On the other hand, Cochran and Harvey defined a container $\F_n(G) =
\tilde G_n$ of $\H_n$ as follows.  Let $\tilde G_0 = \{e\}$, a trivial
group.  Suppose $\tilde G_n$ has been defined as a PTFA group.  Then
$\tilde G_{n+1}$ is defined to be $\tilde G_{n+1} = H_1(G;\K\tilde
G_n) \rtimes \tilde G_n$, where the semidirect product is formed by
viewing $H_1(G;\K\tilde G_n)$ as a $\Z[\tilde G_n]$-module.  Also,
injections $\sigma_{n,G} \colon \H_n(G)=G/G_H\^n \to \tilde
G_n=\F_n(G)$ are defined as follows.  Initially $\sigma_{0,G}$ is the
trivial homomorphism.  Suppose $\sigma_{n,G}$ has been defined.
Consider the composition
\[
\phi_{n,G}\colon G_H\^n/G_H\^{n+1} \xrightarrow{\Phi_{n,G}} H_1(G;\K\H_n(G))
\xrightarrow{f} H_1(G;\K\F_n(G))
\]
where $\Phi_{n,G}$ is the injection induced by the homomorphism used
above to define $G_H\^{n+1}$, and $f$ is induced by $\sigma_{n,G}
\colon \H_n(G) \to \F_n(G)$.  In~\cite{Cochran-Harvey:2004-1} the
followings were shown: there is a derivation $G \to H_1(G;\K\F_n(G))$
which induces~$\phi_{n,G}$.  This derivation, together with $G \to
\H_n(G) \to \F_n(G)$, gives rise to a homomorphism $G \to
H_1(G;\K\F_n(G)) \rtimes \F_n(G) = \F_{n+1}(G)$ with kernel
$G_H\^{n+1}$.  So it induces an injection $\sigma_{n+1,G}\colon
\H_{n+1}(G) \to \F_{n+1}(G)$.  For $n=\omega$, $\tilde G_\omega$
(which is denoted by $\tilde G$ and called the solvable completion of
$G$ in \cite{Cochran-Harvey:2004-1}) is defined to be $\tilde G_\omega
= \varprojlim \tilde G_k$ where $k<\omega$.  For any $n\le \omega$,
$\sigma_n \colon \H_n \to \F_n$ is a container of $\H_n$.  For further
details see~\cite{Cochran-Harvey:2004-1}.

We note that from the definitions it follows that there is a
commutative diagram
\[
\tag*{$(*\mathord{*}*)$}
\begin{diagram}\dgARROWLENGTH=1.8em \dgHORIZPAD=.2em
  \node{\kern 2em 1} \arrow{e}
  \node{G_H\^n/G_H\^{n+1}} \arrow{e} \arrow{s,r}{\phi_{n,G}}
  \node{\H_{n+1}(G)} \arrow{e} \arrow{s,r}{\sigma_{n+1,G}}
  \node{\H_n(G)} \arrow{e} \arrow{s,r}{\sigma_{n,G}}
  \node{1\kern 2em}
  \\
  \node{\kern 2em 1} \arrow{e}
  \node{H_1(G;\K\F_n(G))} \arrow{e}
  \node{\F_{n+1}(G)} \arrow{e}
  \node{\F_n(G)} \arrow{e}
  \node{1\kern 2em}
\end{diagram}
\]
with exact rows.

In the remaining part of this section, we compare the container $\F_n$
with the universal container $\hat\H_n$ of $\H_n$.  More precisely, by
Corollary~\ref{corollary:existence-of-container-of-C-H-quotient},
there is an injective morphism $\hat\sigma_{n,G}$ of our universal
container $\hat\H_n(G) = \hat G / \hat G_H\^n$ to $\F_n(G)=\tilde
G_n$.  Then our question is whether $\hat\sigma_{n,G}$ is an
isomorphism.  The following proposition says that this is closely
related to the structure of certain homology modules of $\hat G$ in an
inductive manner.

\begin{proposition}
  \label{proposition:criterion-for-universality-of-C-H-container}
  For a finitely presented group $G$, $\hat\sigma_{n+1,G}$ is an
  isomorphism if and only if $\hat\sigma_{n,G}$ is an isomorphism and
  the canonical homomorphism
  \[
  H_1(\hat G;\Z\H_n(\hat G)) \to H_1(\hat G;\K\H_n(\hat G))
  \]
  is surjective.
\end{proposition}

\begin{proof}
  Before proving the proposition, we assert that if
  $\hat\sigma_{n,G}\colon \H_n(\hat G) \to \F_n(G)$ is an isomorphism
  then there is a commutative diagram
  \[
  \begin{diagram}\dgHORIZPAD=.2em
    \node{\kern 2em 1} \arrow{e}\dgARROWLENGTH=2em
    \node{\hat G_H\^n/\hat G_H\^{n+1}}
      \arrow{e} \arrow{s,r}{\Phi_{n,\hat G}}
    \node{\H_{n+1}(\hat G)} \arrow{e} \arrow{s,r}{\hat \sigma_{n+1,G}}
    \node{\H_n(\hat G)} \arrow{e} \arrow{s,r}{\hat \sigma_{n,G}}
    \node{1\kern 2em}
    \\
    \node{\kern 2em 1} \arrow{e}
    \node{H_1(\hat G;\K\H_n(\hat G))} \arrow{e}
    \node{\F_{n+1}(G)} \arrow{e}
    \node{\F_n(G)} \arrow{e}
    \node{1\kern 2em}
  \end{diagram}
  \]
  with exact rows. 

  To prove this, we recall that the morphism $\hat\sigma_n$ can be
  described as follows.  Choose a sequence $G = G_0 \to G_1 \to
  \cdots$ of morphisms in $\Omega^\Q$ whose limit is $\hat G$.  Then
  $\H_n(\hat G) \cong \varinjlim \H_n(G_k)$ and
  $\hat\sigma_{n,G}\colon \H_n(\hat G) \to \F_n(G)$ is the limit of
  \[
  \sigma_{n,G_k}\colon \H_n(G_k) \to \F_n(G_k) \cong \F_n(G)
  \]
  as $k\to \infty$.

  By the hypothesis of the assertion, $\H_n(\hat G) \cong \F_n(G)
  \cong \F_n(G_k)$ where the latter isomorphism is induced by $G \to
  G_k$.  So from the diagram $(*\mathord{*}*)$ we obtain a commutative
  diagram%
  {\footnotesize
  \[\hbox to 0mm{\hss$
  \begin{diagram}\dgARROWLENGTH=-6em\dgHORIZPAD=.4em
    \node{\kern 1em 1} \arrow{e}
    \node{\frac{G_H\^n}{G_H\^{n+1}}}
      \arrow[2]{e} \arrow{sse} \arrow[4]{s,l,3}{\phi_{n,G}}
    \node[2]{\H_{n+1}(G)}  \arrow[2]{e} \arrow{sse} \arrow[2]{s,-}
    \node[2]{\H_{n}(G)}  \arrow{e} \arrow{sse} \arrow[2]{s,-}
    \node{1 \kern 1em}
    \\
    \\
    \node[2]{\kern 1em 1} \arrow{e}
    \node{\frac{(G_k)_H\^n}{(G_k)_H\^{n+1}}}
      \arrow[2]{e} \arrow[4]{s,r,3}{\phi_{n,G_k}}
    \node{} \arrow[2]{s,r}{\sigma_{n+1,G}}
    \node{\H_{n+1}(G_k)}  \arrow[2]{e} \arrow[4]{s,r,3}{\sigma_{n+1,G_k}}
    \node{} \arrow[2]{s,r}{\sigma_{n,G}}
    \node{\H_{n}(G_k)}  \arrow{e} \arrow[4]{s,r,3}{\sigma_{n,G_k}}
    \node{1 \kern 1em}
    \\
    \\
    \node{\kern 1em 1} \arrow{e}
    \node{H_1(G;\K\H_n(\hat G))} \arrow{e,-} \arrow{sse}
    \node{} \arrow{e}
    \node{\F_{n+1}(G)}  \arrow{e,-} \arrow{sse,b}{\cong}
    \node{} \arrow{e}
    \node{\F_{n}(G)}  \arrow{e} \arrow{sse,b}{\cong}
    \node{1\kern 1em}
    \\
    \\
    \node[2]{\kern 1em 1} \arrow{e}
    \node{H_1(G_k;\K\H_n(\hat G))} \arrow[2]{e}
    \node[2]{\F_{n+1}(G_k)} \arrow[2]{e}
    \node[2]{\F_n(G_k)} \arrow{e}
    \node{1 \kern 1em}
  \end{diagram}
  $\hss}
  \]
  }%
  with exact rows.  $H_1(G;\K\H_n(\hat G)) \to H_1(G_k;\K\H_n(\hat
  G))$ is an isomorphism since $\F_{n+1}(G) \to \F_{n+1}(G_k)$ and
  $\F_{n}(G) \to \F_{n}(G_k)$ are isomorphisms.  Since $H_1$ commutes
  with limits,
  \[
  H_1(G;\K\H_n(\hat G)) = \varinjlim H_1(G_k;\K\H_n(\hat G)) =
  H_1(\hat G;\K\H_n(\hat G)).
  \]
  Moreover, from the limit of the second row it follows that
  \[
  \varinjlim\, (G_k)_H\^n/(G_k)_H\^{n+1} = \hat G_H\^n / \hat
  G_H\^{n+1}.
  \]
  Also the limit homomorphism
  \[
  \varinjlim \phi_{n,G_k} \colon \hat G_H\^n / \hat G_H\^{n+1} \to
  H_1(\hat G; \K\H_n(\hat G))
  \]
  is equal to the homomorphism $\Phi_{n,\hat G}$.  So the commutative
  diagram in our assertion is obtained by taking the limit.  This
  completes the proof of the assertion.

  Now we prove the proposition.  For the if part, note that it
  suffices to investigate the surjectivity of $\hat\sigma_{n+1,G}$
  since it is always injective.  Since $\hat\sigma_{n,G}$ is an
  isomorphism by the hypothesis, from our assertion it follows that
  $\hat\sigma_{n+1,G}$ is surjective if and only if $\Phi_{n,\hat G}$
  is surjective.  Since $\Phi_{n,\hat G}$ is induced by the
  composition
  \[
  \hat G_H\^n \xrightarrow{p} \hat G_H\^n/[\hat G_H\^n,\hat G_H\^n] =
  H_1(\hat G;\Z\H_n(\hat G)) \xrightarrow{q} H_1(\hat G;\K\H_n(\hat
  G))
  \]
  and $p$ is surjective, $\Phi_{n,\hat G}$ is surjective if and only
  if $q$ is surjective.  This proves the if part.

  For the only if part, note that $\hat\sigma_{n+1,G}$ induces
  $\hat\sigma_{n,G}$ on quotient groups.  (To prove this, one may use
  the argument of the proof of our assertion above; it shows that the
  right square in the diagram of the assertion commutes even without
  assuming that $\hat\sigma_{n,G}$ is an isomorphism.)  Since
  $\hat\sigma_{n+1,G}$ is surjective by the hypothesis, so is
  $\hat\sigma_{n,G}$.  Since $\hat\sigma_{n,G}$ is always injective,
  it is an isomorphism.  So from our assertion it follows that
  $H_1(\hat G;\Z\H_n(\hat G)) \to H_1(\hat G;\K\H_n(\hat G))$ is
  surjective as in the previous paragraph.  This proves the only if
  part.
\end{proof}

Now we investigate inductively whether $\hat\sigma_{n,G}$ is an
isomorphism.

\begin{proposition}
  \label{proposition:C-H-container-for-low-index}
  $\hat\sigma_{n,G}$ is an isomorphism between $\hat\H_n(G)$ and
  $\F_n(G)$ for $n<2$.
\end{proposition}

\begin{proof}
  For $n=0$, $\hat\sigma_{0,G}$ is obviously an
  isomorphism, being a homomorphism between trivial groups.
  
  For $n=1$, we consider
  \[
  \Phi_{0,\hat G}\colon H_1(\hat G;\Z) \to H_1(\hat G;\Q).
  \]
  By the lemma below, which is a special case of
  Lemma~\ref{lemma:divisibility-of-H_1-of-R-closed-group} proved
  later, it follows that $\Phi_{0,\hat G}$ is an isomorphism:
  
  \begin{lemma}
    \label{lemma:divisibility-of-H_1-of-hat-G}
    For any group $G$, $H_1(\hat G;\Z)$ is divisible.
  \end{lemma}

  By
  Proposition~\ref{proposition:criterion-for-universality-of-C-H-container},
  the proof of
  Proposition~\ref{proposition:C-H-container-for-low-index} is
  completed.
\end{proof}

However, for $n=2$, the following result illustrates that $\hat
\sigma_n \colon \hat\H_n \to \F_n$ is not necessarily an isomorphism:

\begin{proposition}
  \label{proposition:non-surjectivity-for-free-groups}
  For any free group $F$ with rank $>1$, $\hat \sigma_{2,F} \colon
  \hat\H_2(F) \to \F_2(F)$ is not surjective.
\end{proposition}

From Proposition \ref{proposition:non-surjectivity-for-free-groups},
it follows that $\hat \sigma_{n,F} \colon \hat\H_n(F) \to \F_n(F)$ is
not surjective for any $n > 1$ (including $n=\omega$), since $\hat
\sigma_{n,F}$ induces a non-surjective homomorphism, namely
$\hat\sigma_{2,F}$, on quotient groups.  Therefore we obtain the
following result:

\begin{theorem}
  \label{theorem:non-universality-of-C-H-container}
  The container $\F_n$ of $\H_n$ is not universal for any $2 \le n \le
  \omega$.
\end{theorem}

\begin{proof}[Proof of
  Proposition~\ref{proposition:non-surjectivity-for-free-groups}]
  We start with a general discussion about an arbitrary finitely
  presented group~$G$.  By
  Lemma~\ref{lemma:divisibility-of-H_1-of-hat-G}, $H_1(\hat G;\Z)$ is
  divisible.  In fact, $H_1(\hat G;\Z)$ is a $\Q$-module (see
  Lemma~\ref{lemma:divisibility-of-H_1-of-R-closed-group}).  So, by
  the definition, $\hat G_H\^1$ is the kernel of the surjection $\hat
  G \to H_1(\hat G;\Q) = H_1(\hat G;\Z)\otimes \Q = H_1(\hat G;\Z)$.
  It follows that $\hat G_H\^1$ is equal to the ordinary commutator
  subgroup $\hat G\^1 = [\hat G,\hat G]$.  Also,
  \[
  \H_1(\hat G)=\hat G/\hat G_H\^1 = H_1(\hat G;\Q)=H_1(G;\Q)=\Q^\mu,
  \]
  where $\mu$ is the first betti number of~$G$.  (The third equality
  is a well-known property of a homology localization.  For
  concreteness we remark that it can be shown by appealing to Theorem
  \ref{theorem:properties-of-closures} (2): there is a sequence $G=G_0
  \to G_1 \to \cdots$ of rationally 2-connected homomorphisms with
  limit $\hat G$.)

  We consider
  \[
  \Psi \colon H_1(\hat G; \Z H_1(\hat G;\Q)) \to H_1(\hat G;
  \K H_1(\hat G;\Q)).
  \]
  Since $H_1(\hat G;\Q)$ is abelian, $\K H_1(\hat G;\Q)$ is the
  ordinary localization $S^{-1} \cdot \Z H_1(\hat G;\Q)$ of the
  commutative ring $\Z H_1(\hat G;\Q)$ where $S=\Z H_1(\hat
  G;\Q)-\{0\}$.  Moreover
  \[
  H_1(\hat G; \K H_1(\hat G;\Q)) = S^{-1} \cdot
  H_1(\hat G; \Z H_1(\hat G;\Q))
  \]
  since $S^{-1} \cdot \Z H_1(\hat G;\Q)$ is a flat $\Z H_1(\hat
  G;\Q)$-module.  Therefore $\Psi$ is surjective if and only if every
  element in $H_1(\hat G; \Z H_1(\hat G;\Q))$ is divisible by any
  element in $S$, that is, for any $u \in H_1(\hat G; \Z H_1(\hat
  G;\Q))$ and $s\in S$, there exists $v \in H_1(\hat G; \Z H_1(\hat
  G;\Q))$ such that $s\cdot v = u$.

  We will show that the divisibility criterion is not satisfied in
  case that $G$ is a free group $F$ of rank $\mu > 1$.  Let $x$ and
  $y$ be two distinct generators of~$F$.  As an abuse of notation, for
  an element $g$ in $F$, we denote the image of $g$ under $F \to \hat
  F$ by~$g$.  (Indeed, it can be shown that $F \to \hat F$ is
  injective for a free group~$F$, although we will not use it.)
  Consider the element
  \[
  u \in \hat F\^1 / [\hat F\^1, \hat F\^1] = H_1(\hat F; \Z H_1(\hat
  F;\Q))
  \]
  which is represented by $xyx^{-1}y^{-1} \in \hat F\^1$, and the
  element
  \[
  s = [x]-1 \in \Z H_1(F;\Q)=\Z H_1(\hat F;\Q)
  \]
  where $[x]$ is the homology class of~$x$.  Suppose that there exists
  $v \in \hat F\^1$ such that $s\cdot v = u$ in $H_1(\hat F; \Z
  H_1(\hat F;\Q))$.  Since the action of $H_1(\hat F;\Q)$ on $H_1(\hat
  F;\Z H_1(\hat F;\Q))$ is given by conjugation, we have
  \[
  xvx^{-1}v^{-1} \equiv u = xyx^{-1}y^{-1} \mod [\hat F\^1, \hat F\^1].
  \]
  From this it follows that $xyx^{-1}y^{-1}$ is in $\hat F_3 = [\hat
  F, [\hat F, \hat F]]$, the third term of the lower central series
  of~$\hat F$.
  
  We will show that a contradiction is derived from this.  Obviously,
  $xyx^{-1}y^{-1}$ is not in $F_3$, say by Hall's basis theorem.  To
  generalize this to $\hat F$, we use the rational version of
  Stallings' theorem: the rational derived series $G_q^\Q$ of a group
  $G$ is defined inductively by
  \[
  G_1^\Q = G, \quad G_{q+1}^\Q = \text{kernel of } G_q^\Q \to
  \frac{G_q^\Q}{[G,G_q^\Q]} \to \frac{G_q^\Q}{[G,G_q^\Q]}
  \mathbin{\mathop{\otimes}_\Z} \Q.
  \]
  Then for any group homomorphism $\pi \to G$ which is rationally
  2-connected, $\pi/\pi_q^\Q \to G/G_q^\Q$ is injective for
  all~$q$~\cite{Stallings:1965-1}.  We also need the following facts:
  obviously $G_q \subset G_q^\Q$ for any group $G$, and for a free
  group $F$, $F_q = F_q^\Q$ since $F_q/F_{q+1}$ is known to be torsion
  free as an abelian group.

  Now, applying the rational version of Stallings' theorem to our $F
  \to \hat F$ which is rationally 2-connected, it follows that
  $xyx^{-1}y^{-1} \in F_3^\Q = F_3$ since $xyx^{-1}y^{-1} \in \hat F_3
  \subset \hat F_3^\Q$.  This is a contradiction.
\end{proof}

\begin{remark}
  In the proof of
  Proposition~\ref{proposition:non-surjectivity-for-free-groups}, we
  considered a particular element $s=[x]-1$ in $S$ to show that the
  divisibility criterion is not satisfied.  Our argument also works
  for any $s$ contained in the kernel of the augmentation homomorphism
  $\Z H_1(F;\Q) \to \Z$.  Such an element $s$ can be used for this
  purpose since it is invertible in the Ore localization; but it is
  not in the Cohn localization.  As mentioned in the introduction,
  this fact motivates a study of a more natural series similar to the
  torsion-free derived series but defined using the Cohn localization
  instead of the Ore localization.  (See also Remark 5.22 of
  \cite{Cochran-Harvey:2004-1}.)  We do not address this issue in
  depth in the present paper.
\end{remark}

\section{Nullhomologous equations and $R$-closures}
\label{section:nullhomologous-equations}

Let $G$ be a group.  We call an element $w$ in the free product
$G*F\langle x_1,\ldots,x_n\rangle$ a \emph{monomial over $G$ in
  $x_1,\ldots,x_n$}, where $F\langle x_1,\ldots,x_n\rangle$ denotes
the free group generated by $x_1,\ldots, x_n$.  Viewing a monomial $w$
as a word in elements of $G$ and $x_1,\ldots,x_n$, we sometimes write
$w=w(x_1,\ldots,x_n)$.

We consider systems of equations over $G$ of the following form:
\[
x_i^e = w_i(x_1,\ldots,x_n), \quad i=1,\ldots, n
\]
where $x_1,\ldots,x_n$ are considered as indeterminates, $e$ is a
nonzero integer, and $w_i(x_1,\ldots,x_n)$ is a monomial over~$G$.  An
$n$-tuple $(g_1,\ldots, g_n)$ of elements in $G$ is called a
\emph{solution} of the system if $x_i = g_i$ satisfies the equations,
that is, $g_i^e$ is equal to $w_i(g_1,\ldots,g_n)$ in $G$ for all~$i$.

Henceforth we fix a subring $R$ of~$\Q$.  We denote by $D_R$ the set
of denominators of reduced fractional expressions of elements in~$R$.
$D_R$~is a multiplicatively closed set.

\begin{definition}
  A system $\{x_i^e = w_i(x_1,\ldots,x_n)\}$ over $G$ is called
  \emph{$R$-nullhomologous} if $e$ is in $D_R$ and each
  $w_i(x_1,\ldots,x_n)$ is sent to the trivial element by the
  canonical projection
  \[
  G*F \to F \to F/[F,F] = H_1(F)
  \]
  where $F=F\langle x_1,\ldots,x_n\rangle$.
\end{definition}

\begin{definition}
  \begin{enumerate}
  \item A group $A$ is called \emph{$R$-closed} if every
    $R$-nullhomologous system over $A$ has a unique solution in~$A$.
  
  \item For a group $G$, an $R$-closed group $\hat G$ equipped with a
    homomorphism $G \to \hat G$ is called an \emph{$R$-closure} of $G$
    if for any homomorphism of $G$ into an $R$-closed group $A$, there
    exists a unique homomorphism $\hat G \to A$ making the following
    diagram commute:
    \[
    \begin{diagram}
      \node{G}\arrow{s}\arrow{e}
      \node{\hat G}\arrow{sw,..}\\
      \node{A}
    \end{diagram}
    \]
    That is, $G\to \hat G$ is the universal (initial) object in the
    category of homomorphisms of $G$ into $R$-closed groups.
  \end{enumerate}
\end{definition}

\begin{remark}
  Although we do not need it in this paper, it can be seen that every
  $R$-nullhomologous system has a unique solution if and only if so
  does every system of the form $\{x_i^e = g_i u_i(x_1,\ldots,x_n)\}$
  where $g_i \in G$, $u_i \in [F,F]$, and $e\in D_R$.  This form is
  more similar to the equations considered in work of
  Farjoun--Orr--Shelah~\cite{Farjoun-Orr-Shelah:1989}.  The only if
  part is clear.  For the if part, suppose an $R$-nullhomologous
  system $\{x_i^e = w_i(x_1,\ldots,x_n)\}$ over $G$ is given.  If the
  variables $x_i$ commuted with elements of $G$ appearing in $w_i$,
  then $w_i$ would be of the form $g_i \cdot u_i$ where $g_i \in G$
  and $u_i\in [F,F]$.  Therefore we can rewrite $w_i$ as $g_i \cdot
  \big(\prod_j [h_{ij}, x_{ij}^{\pm1}]\big) \cdot u_i$ where $h_{ij}
  \in G$.  For each $h_{ij}$, we adjoin to the system an indeterminate
  $y_{ij}$ and an equation $y_{ij} = h_{ij}$.  Replacing each
  occurrence of $h_{ij}$ in the original equation $x_i^e=w_i$ by the
  new indeterminate $y_{ij}$, we obtain a system of the desired form.
  From this the assertion follows.
\end{remark}

The following definition generalizes the notion of invisible subgroups
in Levine's work~\cite{Levine:1989-1}.

\begin{definition}
  A normal subgroup $N$ in $G$ is called \emph{$R$-invisible} if
  \begin{enumerate}
  \item $N$ is normally finitely generated in $G$, and
  \item The order of every element in $N/[G,N]$ is (finite and)
    in~$D_R$.
  \end{enumerate}
\end{definition}

We recall that $N$ is said to be \emph{normally finitely generated} in
$G$ if there exist finitely many elements $a_1,\ldots,a_n$ in $G$
such that $N$ is the smallest normal subgroup containing the~$a_i$.
In this case the $a_i$ are called normal generators of~$N$.  Note that
a normally finitely generated subgroup $N$ in $G$ is $R$-invisible if
and only if $(N/[G,N])\otimes_\Z R =0$.

\begin{remark}
  In his work on homology localizations, Bousfield called a normal
  subgroup $N$ in $G$ \emph{$\pi$-perfect} if $N=[G,N]$
  \cite{Bousfield:1975-1}.  For our purpose, we need to modify it
  regarding the coefficient $R$ and the finiteness assumption as in the
  above definition.  When $R=\Z$, our definition agrees with the
  definition of an \emph{invisible} subgroup due to
  Levine~\cite{Levine:1989-1}.
\end{remark}

In what follows we discuss some useful relationships between
$R$-invisible subgroups and $R$-nullhomologous systems.

\begin{lemma}\label{lemma:invisible-subgps-are-killed-in-closed-gps}
  Suppose $\phi\colon G \to A$ is a homomorphism into an $R$-closed
  group~$A$.  Then every $R$-invisible subgroup $N$ in $G$ is
  contained in the kernel of~$\phi$.
\end{lemma}

\begin{proof}
  Choose normal generators $a_1,\ldots,a_n$ of $N$.  Since the order
  of $a_i$ is in $D_R$, there exist an element $e\in D_R$ such that
  $a_i^e \in [G,N]$ for all~$i$.  So we can write $a_i^e$ as a product
  of commutators $[b_{ij},c_{ij}]$ where $b_{ij}\in N$, $c_{ij} \in
  G$.  Furthermore $b_{ij}$ can be written as a product of conjugates
  of the~$a_k$.  Replacing each occurrence of $a_k$ in this expression
  of $b_{ij}$ by an indeterminate $x_k$ and plugging the result into
  the above expression of $a_i^e$, we obtain a word
  $w_i(x_1,\ldots,x_n)$ in $G*F$, $F$ is the free group generated by
  the $x_i$, such that the system
  \[
  x_i^e = w_i(x_1,\ldots,x_n), \qquad i=1,\ldots,n
  \]
  has two sets of solutions, $\{x_i = 1\}$ and $\{x_i = a_i\}$.  It is
  easily seen that this system is $R$-nullhomologous.

  Denote by $w_i^\phi$ the image of $w_i$ under $G*F \to A*F$.  It
  gives rise to an $R$-nullhomologous system
  \[
  x_i^e = w_i^\phi(x_1,\ldots,x_n), \qquad i=1,\ldots,n
  \]
  over $A$, which has two solution sets $\{x_i = 1\}$ and $\{x_i =
  \phi(a_i)\}$.  By the uniqueness of a solution over $A$, it follows
  that $\phi(a_i)=1$.  Thus $\phi(N)$ is trivial.
\end{proof}

\begin{lemma}\label{lemma:product-of-R-invisible-subgps}
  If $N_1$ and $N_2$ are $R$-invisible subgroups in $G$, then
  $N=N_1N_2$ is also $R$-invisible in~$G$.
\end{lemma}

\begin{proof}
  Obviously $N$ is normally finitely generated.  For any $n\in N$,
  write $n=n_1n_2$ where $n_i \in N_i$.  Since $N_i$ is $R$-invisible,
  there exist $e\in D_R$ such that $n_i^e \in [G,N_i]$.  Then
  \[
  (n_1n_2)^e \equiv n_1^e n_2^e \equiv 1 \mod [G,N].
  \]
  This shows that the order of $n[G,N]$ in $N/[G,N]$ is (a divisor
  of)~$e$.
\end{proof}

\begin{lemma}\label{lemma:no-nontrivial-invisible-subgroup}
  Suppose $G$ is a group and $N$ is the union of all $R$-invisible
  subgroups in~$G$.  Then $G/N$ has no nontrivial $R$-invisible
  subgroup.
\end{lemma}

\begin{proof}
  First of all, $N$ is a normal subgroup by the previous lemma.
  Suppose $H/N$ is $R$-invisible in $G/N$ for some $H\subset G$.
  Choose a finite normal generator set $\{h_iN \}$ of $H/N$, $h_i \in
  H$.  It suffices to show that $h_i \in N$ for each~$i$.  For some
  $e\in D_R$, $h_i^e \in [H/N,G/N]$.  Therefore $h_i^e \in n_i[H,G]$
  for some $n_i \in N$.  By the previous lemma, there exists an
  $R$-invisible subgroup $K$ in $G$ such that $n_i \in K$ for all~$i$.
  Then the normal subgroup $K_1$ generated by $K$ and the $h_i$ is
  $R$-invisible in~$G$.  It follows that $h_i \in K_1 \subset N$.
\end{proof}

\begin{lemma}\label{lemma:uniqueness-of-solution}
  If $G$ has no nontrivial $R$-invisible subgroup, then any
  $R$-nullhomologous system over $G$ has at most one solution.
\end{lemma}

\begin{proof}
  Suppose $S=\{x_i^e = w_i(x_1,\ldots,x_n)\}$ is $R$-nullhomologous
  and $\{x_i=a_i\}$ and $\{x_i=b_i\}$ are solutions of $S$ over~$G$.
  Let $N$ be the normal subgroup in $G$ generated by the
  $a_i^{\mathstrut} b_i^{-1}$, and let
  \[
  u_i(x_1,\ldots,x_n)=w_i(x_1b_1,\ldots,x_nb_n)b_i^{-e}.
  \]
  Since $u_i(1,\ldots,1) = w_i(b_1,\ldots,b_n)b_i^{-e} = 1$,
  we can write $u_i$ as a word of the form
  \begin{align*}
    u_i &= \prod_j g_{ij}^{\mathstrut} x_{ij}^{\pm1} g_{ij}^{-1} \\
    &= \Big( \prod_j (x_{i1}^{\pm1} \cdots x_{i,j-1}^{\pm1})
    [g_{ij}^{\mathstrut}, x_{ij}^{\pm1}] (x_{i1}^{\pm1} \cdots
    x_{i,j-1}^{\pm1})^{-1} \Big) \prod_j x_{ij}^{\pm1}
  \end{align*}
  where $g_{ij} \in G$ and $x_{ij} = x_{k_{ij}}$ for some~$k_{ij}$.
  Furthermore, each $u_i$ is killed by $G*F \to F/[F,F]$ where
  $F=F\langle x_1,\ldots,x_n \rangle$, since so is~$w_i$.  It follows
  that $\prod_j x_{ij}^{\pm1}$, the last term of the above expression,
  is contained in $[F,F]$.  Therefore
  \begin{align*}
    a_i^{e} b_i^{-e} &= w_i(a_1,\ldots,a_n)b_i^{-e} \\
    &= u_i(a_1^{\mathstrut} b_1^{-1},\ldots, a_n^{\mathstrut} b_n^{-1}) \in [G,N].
  \end{align*}
  Now we have
  \begin{align*}
    (a_i^{\mathstrut} b_i^{-1})^e &\equiv (a_i^{\mathstrut} b_i^{-1})^e \cdot a_i^{-e} b_i^{e} \\
    & \equiv (a_i^{\mathstrut}b_i^{-1})^{e-1} \cdot a_i^{\mathstrut}b_i^{-1} \cdot a_i^{-e} b_i^{e} \\
    & \equiv (a_i^{\mathstrut}b_i^{-1})^{e-1} \cdot a_i^{-e} \cdot a_ib_i^{-1} \cdot b_i^{e} \\
    & \equiv (a_i^{\mathstrut}b_i^{-1})^{e-1} \cdot a_i^{-e+1} b_i^{e-1} \equiv
    \cdots \equiv 1 \mod [G,N]
  \end{align*}
  This shows that $N$ is $R$-invisible.  By the hypothesis, $N$ is
  trivial and $a_i=b_i$.
\end{proof}

As an immediate consequence of the definition, we have the following
divisibility result:

\begin{lemma}
  \label{lemma:divisibility-of-H_1-of-R-closed-group}
  If $A$ is $R$-closed, then $H_1(A;\Z)$ is an $R$-module.
\end{lemma}

\begin{proof}
  Let $g$ be an element in $A$, and let $e$ be an element in~$D_R$.
  Consider the equation $x^e = g$.  Since it is $R$-nullhomologous,
  there is a solution $x=h$ in~$A$.  It follows that the homology
  class of $g$ is divisible by~$e$.
\end{proof}

\section{$R$-closures and localizations}
\label{section:localization-wrt-2-connected-morphisms}

We begin this section by recalling the definition of a localization.
In general, we think of a category $\mathcal{C}$ and a class of
morphisms~$\Omega$ in $\mathcal{C}$.

\begin{definition}
  \begin{enumerate}
  \item An object $A$ in $\mathcal{C}$ is called \emph{local} with
    respect to $\Omega$ if for any morphism $\pi \to G$ in $\Omega$
    and any morphism $\pi \to A$, there exists a unique morphism $G
    \to A$ making
    \[
    \begin{diagram}
      \node{\pi}\arrow{e}\arrow{s}
      \node{G} \arrow{sw,..} \\
      \node{A}
    \end{diagram}
    \]
    commute.

  \item A \emph{localization} with respect to $\Omega$ is a pair
    $(E,p)$ of a functor $E\colon \mathcal{C} \to \mathcal{C}$ and a
    natural transformation $p\colon \id_\mathcal{C} \to E$ (that is,
    each object $G$ is equipped with a morphism $p_G\colon G \to
    E(G)$) such that for any morphism $G\to A$ into a local object $A$,
    there is a unique morphism $E(G) \to A$ making
    \[
    \begin{diagram}
      \node{G}\arrow{e}\arrow{s}
      \node{E(G)} \arrow{sw,..} \\
      \node{A}
    \end{diagram}
    \]
    commute.
  \end{enumerate}
\end{definition}

Of course our main interest is the localization of groups with respect
to the class $\Omega^R$; recall that $\Omega^R$ is the class of group
homomorphisms $\phi\colon \pi \to G$ where $\pi$ is finitely
generated, $G$ is finitely presented, and $\phi$ is 2-connected on
$R$-homology.  From now on local groups and localizations are always
with respect to our~$\Omega^R$.

\begin{theorem}\label{theorem:closed-iff-local}
  A group $A$ is $R$-closed if and only if $A$ is local with respect
  to~$\Omega^R$.
\end{theorem}

\begin{proof}
  Suppose that $A$ is $R$-closed.  To show that $A$ is local, suppose
  that a morphism $\alpha\colon \pi \to G$ in $\Omega^R$ and a
  morphism $\phi\colon \pi \to A$ are given.  We will show that there
  exists a unique morphism $\varphi\colon G \to A$ that $\phi\colon
  \pi\to A$ factors through.

  Choose generators $h_1,\ldots,h_n$ of~$G$.  Since $\alpha$ induces
  an isomorphism 
  \[
  H_1(\pi;R)=\frac{\pi}{[\pi,\pi]}
  \mathbin{\mathop{\otimes}\limits_\Z} R \to \frac{G}{[G,G]}
  \mathbin{\mathop{\otimes}\limits_\Z} R,
  \]
  it follows that for any element $x$ in $G/[G,G]$, $x^r$ is contained
  in the image of $\pi/[\pi,\pi]$ for some $r\in D_R$.  Therefore
  there exists $e\in D_R$ such that each $h_i^e$ can be written as
  \[
  h_i^e = \alpha(g_i) \cdot \prod_j [u_{ij}, v_{ij}]
  \]
  where $g_i\in \pi$ and $u_{ij}=u_{ij}(h_1,\ldots,h_n)$,
  $v_{ij}=v_{ij}(h_1,\ldots,h_n)$ are words in $h_1,\ldots,h_n$.

  Consider the system $S$ of equations $x_i^e = w_i(x_1,\ldots,x_n)$
  where the element $w_i$ in $G*F\langle x_1,\ldots,x_n\rangle$ is
  given by
  \[
  w_i(x_1,\ldots,x_n) = g_i \cdot \prod_j
  [u_{ij}(x_1,\ldots,x_n),v_{ij}(x_1,\ldots,x_n)].
  \]
  Then $S$ is $R$-nullhomologous.

  We associated to the system $S$ a new group $\pi_S$ obtained by
  ``adding'' to $\pi$ a solution $\{z_i\}$ to $S$; formally, it is
  defined to be a amalgamated product of $\pi$ and $F\langle
  z_1,\ldots,z_n\rangle$:
  \[
  \pi_S = \langle \pi, z_1,\ldots,z_n \mid z_i^e=w_i(z_1,\ldots,z_n),
 \, i=1,\ldots,n\rangle.
  \]
  Note that, since $e\in D_R$, the canonical homomorphism $\pi \to
  \pi_S$ is 2-connected on $R$-homology; it can be seen by computing
  $H_1(\pi_S;R)$ and $H_2(\pi_S;R)$ using the complex obtained by
  attaching to $K(\pi,1)$ 1-cells and 2-cells corresponding the new
  generators and relations.

  $\pi \to G$ and $z_i \to h_i \in G$ induce a surjection $\beta
  \colon \pi_S \to G$.  Since $A$ is $R$-closed, there is a unique
  solution $\{x_i=a_i\}$ of the $R$-nullhomologous system $S^\phi =
  \{x_i^e = w_i^\phi(x_1,\ldots,x_n)\}$ over~$A$, which is the image
  of $S$ under~$\phi$.  $\phi\colon \pi \to A$ and $z_i \to a_i \in A$
  induce a homomorphism $\phi_S\colon \pi_S \to A$.
  \[
    \begin{diagram}
    \node{\pi}\arrow{e}\arrow{s,l}{\phi}
    \node{\pi_S} \arrow{sw,l}{\phi_S} \arrow{e,t}{\beta}
    \node{G} \arrow{sww,b,..}{\varphi} \\
    \node{A}
  \end{diagram}
  \]

  We will show that $\phi_S$ factors through a homomorphism
  $\varphi\colon G \to A$.  Let $N$ be the kernel of~$\beta$.
  Applying the Stallings exact sequence~\cite{Stallings:1965-1} to
  $\beta$ which is surjective, we obtain a long exact sequence
  \[
  H_2(\pi_S;R) \to H_2(G;R) \to \frac{N}{[\pi_S,N]}\otimes_\Z R \to
  H_1(\pi_S;R) \to H_1(G;R).
  \]
  Since $\alpha \colon \pi \to G$ and $\pi \to \pi_S$ are 2-connected
  on $R$-homology, so is~$\beta$.  Thus $(N/[\pi_S,N])\otimes R=0$.
  Since $\pi_S$ is finitely generated and $G$ is finitely presented,
  $N$ is finitely normally generated in~$\pi_S$.  This shows that $N$
  is $R$-invisible in~$\pi_S$.  By
  Lemma~\ref{lemma:invisible-subgps-are-killed-in-closed-gps},
  $\phi_S(N)$ is trivial.  It follows that $\phi_S$ induces a
  homomorphism $\varphi\colon G\to A$ as desired.

  If another homomorphism $\varphi' \colon G \to A$ satisfies $\phi =
  \varphi'\circ \alpha$, then $\{x_i = \varphi'(h_i)\}$ is a solution of
  the system $S^\phi$.  By the uniqueness of a solution, we have
  $\varphi(h_i) = a_i = \varphi'(h_i)$, that is, $\varphi = \varphi'$.
  This completes the proof of the only if part.

  For the converse, suppose that $A$ is local, and an
  $R$-nullhomologous system $S=\{x_i^e = w_i(x_1,\ldots,x_n)\}$ over
  $A$ is given.  There are finitely many elements in $A$ which appears
  in the words $w_i$.  Let $G$ be the free group generated by (symbols
  associated to) these elements, and $\phi \colon G \to A$ be the
  natural homomorphism.  Lifting $S$, we obtain a system $S' = \{x_i^e
  = w_i(x_1,\ldots,x_n)\}$ over $G$ which is sent to $S$ by~$\phi$.
  Consider the group $G_{S'}$ obtained by ``adding a solution
  $\{z_i\}$ to the system $S'$'' as before.  Then $G_{S'}$ is finitely
  presented and the canonical homomorphism $\alpha \colon G \to
  G_{S'}$ is 2-connected.  Since $A$ is local, there exists a unique
  homomorphism $\varphi\colon G_{S'} \to A$ making
  \[
  \begin{diagram}
    \node{G}\arrow{e,t}{\alpha} \arrow{s,l}{\phi}
    \node{G_{S'}} \arrow{sw,r,..}{\varphi} \\
    \node{A}
  \end{diagram}
  \]
  commute.  Now $\{\varphi(z_i)\}$ is a solution of $S$ over $A$.
  
  If there is another solution $\{x_i = a_i\}$ of $S$, then $x_i \to
  a_i$ gives rise to another homomorphism $\varphi'\colon G_{S'} \to
  A$ making the above diagram commute.  By the uniqueness of
  $\varphi$, $\varphi' =\varphi$, and therefore, $a_i = \varphi(z_i)$.
  This completes the proof.
\end{proof}

\section{Existence of $R$-closures}
\label{section:existence-of-closures}

\begin{theorem}\label{theorem:existence-of-closure}
  For any subring $R$ of $\Q$ and any group $G$, there is an
  $R$-closure $G \to \hat G$.
\end{theorem}

\begin{proof} Basically the construction consists of two parts: adjoin
  solutions repeatedly so that every system has at least one solution
  eventually, and take a quotient of the resulting group to identify
  different solutions if any.

  This idea is formalized as follows.  We construct a sequence $G_0,
  G_1,\ldots$ of groups inductively.  Let $G_0=G$.  Suppose $G_k$ has
  been defined.  Let $\mathcal{S}_k$ be the set of all
  $R$-nullhomologous systems over~$G_k$.  We associate a symbol $z_i$
  to an indeterminate $x_i$ of a system in~$\mathcal{S}_k$, and let
  $F_k$ be the free group generated by all the symbols $z_i$.  Let
  $G_{k+1} = G_k*F_k$ modulo the relations corresponding the systems
  in~$\mathcal{S}_k$, that is, each equation
  $x_i^e=w_i(x_1,\ldots,x_n)$ gives rise to a defining relation
  $z_i^e=w_i(z_1,\ldots,z_n)$ of~$G_{k+1}$.  Let $\bar G=\varinjlim
  G_k$.  Let $N$ be the union of all $R$-invisible subgroups in $\bar
  G$, and finally let $\hat G = \bar G/N$.

  We will show that the canonical homomorphism $\Phi\colon G \to \hat
  G$ is an $R$-closure of~$G$.  First we claim that $\hat G$ is
  $R$-closed.  For the existence of a solution, suppose that
  $S=\{x_i^e = w_i(x_1,\ldots,x_n)\}$ be an $R$-nullhomologous system
  over~$\hat G$.  Since the $w_i$ involve only finitely many elements
  of $\hat G$, $S$ lifts to an $R$-nullhomologous system over some
  $G_k$, that is, each $w_i$ is the image of an element $w_i'$ in $G_k
  * F$ where $F=F\langle x_1,\ldots,x_n\rangle$.  Sending $w_i'$ via
  $G_k *F \to G_{k+1}*F$, we obtain an $R$-nullhomologous system over
  $G_{k+1}$ which has a solution $\{x_i = a_i\}$ by our construction
  of~$G_{k+1}$; recall that $G_{k+1}$ is obtained by adjoining
  solutions of all $R$-nullhomologous systems over~$G_k$.  Obviously
  the image of the $a_i$ in $\hat G$ is a solution of the given
  system~$S$.  On the other hand, $\hat G$ has no nontrivial
  $R$-invisible subgroup by
  Lemma~\ref{lemma:no-nontrivial-invisible-subgroup}.  From
  Lemma~\ref{lemma:uniqueness-of-solution}, the uniqueness of a
  solution follows.  This proves the claim.

  Now it remains to show that $\Phi\colon G\to \hat G$ is a universal
  (initial) object.  Suppose that $\phi\colon G \to A$ is a
  homomorphism of $G$ into an $R$-closed group~$A$.  Since $G_1$ is
  obtained from $G=G_0$ by adjoining solutions of $R$-nullhomologous
  systems, there exists a unique homomorphism $\phi_1\colon G_1 \to A$
  making the below diagram commute; $\phi_1$ is defined by sending new
  generators $z_i$ of $G_1$ associated to a system $S=\{x_i =
  w_i(x_1,\ldots,x_n)\}$ over $G_0$ to the solution of $S^\phi$ over
  $A$, and the uniqueness of $\phi_1$ follows from the uniqueness of a
  solution over~$A$.
  \[
  \begin{diagram}
    \node{G} \arrow{s,l}{\phi}\arrow{e}
    \node{G_1} \arrow{sw,l,..}{\phi_1}\arrow{e}
    \node{G_2} \arrow{sww,r,..}{\phi_2}\arrow{e}
    \node{\cdots} \arrow{e}
    \node{\bar G = \varinjlim G_k} \\
    \node{A}
  \end{diagram}
  \]
  Repeating the same argument, we can inductively construct a sequence
  of homomorphisms $\phi_k\colon G_k \to A$ which make the diagram
  commute.  Passing to the limit, the $\phi_k$ induce a homomorphism
  $\bar\phi\colon \bar G \to A$.  By
  Lemma~\ref{lemma:invisible-subgps-are-killed-in-closed-gps},
  $\varphi'$ kills each $R$-invisible subgroup, and so $N$ is
  contained in the kernel of~$\bar\phi$.  It follows that $\bar\phi$
  gives rise to a homomorphism $\varphi\colon \hat G = \bar G/N \to A$
  such that $\varphi \circ \Phi = \phi$.

  Suppose another homomorphism $\varphi'\colon \hat G \to A$ satisfies
  $\varphi' \circ \Phi = \phi$.  Consider the composition
  $\phi_k'\colon G_k \to \hat G \xrightarrow{\varphi'} A$.  Then the
  $\phi_k$ make the above diagram commute as well as the $\phi_k$.
  From the uniqueness of a solution over $A$, it follows that $\phi_k =
  \phi_k'$ for every~$k$.  Passing to the limit, it follows that
  $\varphi=\varphi'$.  This completes the proof.
\end{proof}

{\long\def\ignoreme{

\begin{remark}
  Compared with the proofs of the existence of algebraic closures for
  the Vogel--Levine localization and the Bousfield localization, the
  most significant difference is that we need to iterate the process
  of adjoining solutions; recall that we have constructed groups $G_k$
  and showed that the algebraic closure is a quotient of the limit.
  In the construction of Levine~\cite{Levine:1989-1} and
  Farjoun--Orr--Shelah~\cite{Farjoun-Orr-Shelah:1989}, it is
  sufficient to consider $G_1$ only, that is, the algebraic closure of
  $G$ is a quotient of the group obtained by adjoining solutions of
  all systems \emph{over $G$}.  The main reason of this sophistication
  in our case is that the exponent $e$ can be greater than~$1$.  In
  fact, if $e$ were always 1, then we could reduce a nullhomologous
  system over $G_{k+1}$ into a system over $G_{k}$ by using a sort of
  back-substitution argument as in \cite{Levine:1989-1}
  and~\cite{Farjoun-Orr-Shelah:1989}, so that every element in the
  limit $\varinjlim G_k$ was (a part of) a solution of a system
  over~$G$. In particular, if the coefficient ring $R$ is $\Z$, then
  our closure has the following property: any element in $\hat G$ is
  (a part of) a solution of a nullhomologous system over~$G$.
\end{remark}

}}

\begin{remark}
  There is an alternative construction of an $R$-closure: let $G_0=G$
  as before, and assuming $G_k$ has been defined, consider the group
  obtained from $G_k$ by adjoining solutions of all $R$-nullhomologous
  systems, and let $G_{k+1}$ be its quotient by the union of all
  $R$-invisible subgroups.  Then it can be proved that $\varinjlim
  G_k$ is an $R$-closure of~$G$.
\end{remark}

\begin{corollary}\label{corollary:existence-of-localization}
  There is a localization $(E,p)$ with respect to $\Omega^R$.
\end{corollary}

\begin{proof}
  For each group $G$, define $E(G)=\hat G$ and $p_G \colon G \to \hat
  G$ to be the homomorphism $\Phi$ constructed above.  For any
  homomorphism $\phi\colon \pi \to G$, the composition $\pi\to G \to
  E(G)=\hat G$ is a homomorphism of $\pi$ into an $R$-closed group.
  By the universal property of the $R$-closure $\pi \to E(\pi)=\hat
  \pi$, there is a unique homomorphism $E(\pi) \to E(G)$ that the
  composition factors through; we denote this homomorphism by
  $E(\phi)\colon E(\pi) \to E(G)$.  One can check that $E$ is a
  functor and $p$ is a natural transformation by using the universal
  property of $R$-closures in a straightforward way.  Finally, since
  $G \to E(G)$ is an $R$-closure, it is initial among homomorphisms of
  $G$ into local groups, by Theorem~\ref{theorem:closed-iff-local}.
\end{proof}

We sometimes denote the induced homomorphism $E(\phi)$ by~$\hat\phi$.

From the universal property of the $R$-closure, a natural isomorphism
theorem follows:

\begin{proposition}
  \label{proposition:isomorphism-induced-by-R-closure}
  If $\alpha\colon\pi \to G$ is in $\Omega^R$, then the induced
  homomorphism $\hat \alpha \colon \hat \pi \to \hat G$ is an
  isomorphism.
\end{proposition}

\begin{proof}
  We consider the following diagram:
  \[
  \begin{diagram}
    \node{\pi} \arrow{e,t}{\alpha} \arrow{s,l}{p_\pi}
    \node{G} \arrow{s,r}{p_G} \arrow{sw,t,..}{\phi} \\
    \node{\hat \pi} \arrow{e,b}{\hat\alpha}
    \node{\hat G}
  \end{diagram}
  \]
  Since $\alpha$ is in $\Omega^R$ and $\hat\pi$ is local, there is
  $\phi\colon G\to \hat \pi$ such that $\phi\alpha=p_\pi$.  Since
  $\hat\pi$ is local and $\hat G$ is the localization, there is
  $\psi\colon \hat G \to \hat\pi$ such that $\psi p_G = \phi$.  We
  will show that $\psi$ is an inverse of $\hat\alpha$.  Observe that
  $(\hat\alpha \psi p_G)\alpha = p_G\alpha$.  Since $\hat G$ is local
  and $\alpha\in \Omega^R$, we have $\hat\alpha \psi p_G = p_G$.  It
  follows that $\hat\alpha\psi=\mathrm{id}_{\hat G}$ by the uniqueness
  of a map $\beta \colon \hat G \to \hat G$ such that $\beta p_G =
  p_G$.  On the other hand, since $\psi \hat\alpha p_\pi = \psi p_G
  \alpha = \phi \alpha = p_\pi$, $\psi \hat\alpha = \mathrm{id}_{\hat
    \pi}$ by a similar uniqueness argument.
\end{proof}

\begin{remark}
  \begin{enumerate}
  \item For subrings $R\subset S\subset \Q$, temporarily denote by
    $\hat G^R$ and $\hat G^S$ the $R$- and $S$-closures of a group
    $G$, respectively.  Then, since $\Omega^R \subset \Omega^S$, there
    is a natural transformation $\hat G^R \to \hat G^S$ making the
    following diagram commute:
    \[
    \begin{diagram}
      \node{G} \arrow[2]{e} \arrow{se}
      \node[2]{\hat G^S} \\
      \node[2]{\hat G^R} \arrow{ne}
    \end{diagram}
    \]
  \item Similar conclusion holds for Levine's algebraic closure
    defined in~\cite{Levine:1989-1}, which is a localization with
    respect to the collection $\Omega^{Levine}$ of group homomorphisms
    $\alpha\colon \pi\to G$ such that $\pi$, $G$ are finitely
    presented, $\alpha$ is (integrally) 2-connected, and $G$ is
    normally generated by the image of~$\alpha$.  Namely, denoting
    Levine's algebraic closure by $\hat G^{Levine}$, there is a
    natural transformation $\hat G^{Levine} \to \hat G^\Z$ making a
    similar diagram commute, since $\Omega^{Levine} \subset
    \Omega^\Z$.
  \end{enumerate}
\end{remark}

We finish this section with some results on the $R$-closure of a
finitely presented group.

\begin{proposition}\label{theorem:closure-as-limit}
  If $G$ is finitely presented, then there is a sequence
  \[
  G=P_0 \to P_1 \to P_2 \to \cdots 
  \]
  of homomorphisms in $\Omega^R$ (in particular, each $P_k$ is
  finitely presented) such that the limit homomorphism $G = P_0 \to
  \varinjlim P_k$ is an $R$-closure, that is, $\hat G\cong \varinjlim
  P_k$.
\end{proposition}

\begin{proof}
  We use the notations of the proof of
  Theorem~\ref{theorem:existence-of-closure}.  Recall that
  $\mathcal{S}_k$ is the set of all $R$-nullhomologous systems over
  $G_k$.  We denote $\mathcal{S}=\bigcup \mathcal{S}_k$.  Then an
  element $S=\{x_i = w_i\}$ in $\mathcal{S}$ is a system over $G_p$
  for some~$p$.  By our construction of $G_p$, each $w_i$ can be
  viewed as a word in indeterminates $x_1,x_2,\ldots$ and solutions of
  other systems over some $G_q$ where $q<p$.

  Since $G_0 = G$ is countable, an induction shows that $G_k$ and
  $\mathcal{S}_k$ are always countable.  Thus the union $\mathcal{S}$
  of all $\mathcal{S}_k$ is countable.  From this it can be seen that
  we can enumerate elements of $\mathcal{S}$ as a sequence $T_1,
  T_2,\ldots$ of systems which satisfies the following: suppose the
  system $T_k = \{x_i=w_i\}$ is over~$G_p$.  Then each $w_i$ involves
  only elements in $G_p$ that can be expressed as a product of
  solutions (and their inverses) of other systems $T_q$ such that
  $q<k$.  In other words, the system $T_{k}$ can be lifted to a system
  over the group
  \[
  Q_{k-1} = (\cdots((G_{T_1})_{T_2})\cdots)_{T_{k-1}}.
  \]
  (Recall our notation that $Q_{k-1}$ is the group obtained from $G$
  by adjoining the solutions of the systems $T_1, T_2, \ldots, T_k$.)
  Note that there is a canonical map $Q_k \to Q_{k+1}$; it is in
  $\Omega^R$ as we discussed before.  Furthermore, it is obvious that
  $\varinjlim Q_k \cong \bar G=\varinjlim G_k$.

  We claim that for each $k$ there is an $R$-invisible subgroup $N_k$
  in $Q_k$ such that
  \begin{enumerate}
  \item $Q_k \to Q_{k+1}$ sends $N_k$ into $N_{k+1}$, and
  \item for any $R$-invisible subgroup $K$ in $\bar G$, there is $k$
    such that $K$ is contained in the normal subgroup of $G$ generated
    by the image of $N_k$ under $Q_k \to \bar G$.
  \end{enumerate}
  For, since $\bar G$ is countable, we can arrange all $R$-invisible
  subgroups in $\bar G$ as a sequence, and by appealing to
  Lemma~\ref{lemma:product-of-R-invisible-subgps}, we can produce an
  increasing sequence $L_1 \subset L_2 \subset \cdots $ of
  $R$-invisible subgroups in $\bar G$ such that every $R$-invisible
  subgroup in $\bar G$ is contained in some~$L_k$.  Since each $L_k$
  $R$-invisible, it is normally generated by finitely many elements
  $b_i$ which satisfy equations of the following form:
  \[
  b_i^e = \prod_j h_j[b_{i_j},g_j]h_j^{-1} \quad \text{where }h_j, g_j
  \in \bar G,\, e\in D_R.
  \]
  All the $b_i,h_j,g_j \in \bar G$ which appear in this expression can
  be lifted to $Q_n$ for some $n$ in such a way that the equations can
  also be lifted, that is, replacing these elements in the equation by
  the lifts, we again obtain an equality in~$Q_n$.  Choosing a
  subsequence of $\{Q_k\}$, we may assume that $n=k$.  Then the lifts
  of $b_i$ in $Q_k$ generate an $R$-invisible subgroup $N_k'$
  in~$Q_k$.  Let $N_k$ be the union of the image of $N_\ell'$ under
  $Q_\ell \to Q_k$ where $\ell$ runs over $1,\ldots,k$.  Now the claim
  follows.
  
  Let $P_k = Q_k/N_k$.  By our claim (1), $Q_k \to Q_{k+1}$ induces
  $P_k \to P_{k+1}$.  We will show that these groups $P_k$ has the
  desired properties.  First, $P_k \to P_{k+1}$ is in $\Omega^R$ since
  $Q_k \to Q_{k+1}$ is 2-connected and so is $Q_k \to P_k$ by the
  Stallings exact sequence.

  So it remains to show that $\varinjlim P_k \cong \hat G$.  For each
  $k$ the composition 
  \[
  j_k \colon Q_k \to \bar G \to \hat G
  \]
  gives rise to $P_k \to \hat G$ since it kills the $R$-invisible
  subgroup~$N_k$.  These morphisms induces $\varinjlim P_k \to \hat
  G$, which is obviously surjective.  To show that it is injective,
  suppose that $j_k$ sends an element $x\in Q_k$ into an $R$-invisible
  subgroup $K$ in~$\bar G$.  We may assume that $K$ is contained in
  the normal closure of $j_k(N_k)$ by our claim (2) above.  Then
  $j_k(x)$ is of the form $\prod_j g_j y_j g_j^{-1}$, $g_j\in \bar G$,
  $y_j \in j_k(N_k)$.  By choosing a sufficiently large $k$, we may
  assume that every $g_j$ is in the image of~$j_k$.  So $j_k(x) \in
  j_k(N_k)$.  By choosing a larger $k$, we may assume that $x \in
  N_k$.  Now it follows that $x$ is sent to the identity in $P_k$.  It
  proves that $\varinjlim P_k \to \hat G$ is injective.
\end{proof}

\begin{proposition}
  Suppose that $\pi_0 \to \pi_1 \to \pi_2 \to \cdots$ and $G_0 \to G_1
  \to G_2 \to \cdots$ are sequences of homomorphisms in $\Omega^R$
  such that $\pi_0 \to \varinjlim \pi_k$ and $G_0 \to \varinjlim G_k$
  are the $R$-closures.  Then for any $\phi \colon \pi_0 \to G_0$,
  there exist a sequence $G_0=P_0 \to P_1 \to P_2 \to \cdots$ of
  homomorphisms in $\Omega^R$ which fits into the commutative diagram
  \[
  \begin{diagram}\dgARROWLENGTH=1em
    \node{\pi_0} \arrow{e} \arrow{s,l}{\phi}
    \node{\pi_1} \arrow{e} \arrow{s}
    \node{\cdots}  \arrow{e}
    \node{\pi_k} \arrow{e} \arrow{s}
    \node{\pi_{k+1}} \arrow{e} \arrow{s}
    \node{\cdots}  \arrow{e}
    \node{\varinjlim \pi_k \hbox to 0mm{ $= \hat \pi$\hss}} \arrow{s}\\
    \node{P_0} \arrow{e}
    \node{P_1} \arrow{e}
    \node{\cdots}  \arrow{e}
    \node{P_k} \arrow{e}
    \node{P_{k+1}} \arrow{e}
    \node{\cdots}  \arrow{e}
    \node{\varinjlim P_k} \\
    \node{G_0} \arrow{e} \arrow{n,=}
    \node{G_1} \arrow{e} \arrow{n}
    \node{\cdots}  \arrow{e}
    \node{G_k} \arrow{e} \arrow{n}
    \node{G_{k+1}} \arrow{e} \arrow{n}
    \node{\cdots}  \arrow{e}
    \node{\varinjlim G_k \hbox to 0mm{ $= \hat G$\hss}} \arrow{n}
  \end{diagram}\hphantom{= \hat G}
  \]
  in such a way that each $G_k \to P_k$ is in $\Omega^R$ and $\hat G
  \to \varinjlim P_k$ is an isomorphism, that is, $G_0 = P_0 \to
  \varinjlim P_k$ is an $R$-closure.
\end{proposition}

\begin{proof}
  We will construct inductively a sequence $n_0 < n_1 < \cdots$ such
  that $P_k = G_{n_k}$ has the desired property.  Let $n_0=0$.  As our
  induction hypothesis, suppose that $n_0,\ldots,n_k$ have been chosen
  in such a way that the following diagram commutes, where $\hat \phi$
  is the map induced by~$\phi$.  (Note that $n_k \ge k$
  automatically.)
  \[
  \begin{diagram}\dgARROWLENGTH=1.6em \dgHORIZPAD=.4em
    \node{\pi_0} \arrow{e} \arrow{s,l}{\phi}
    \node{\cdots} \arrow{e}
    \node{\pi_k} \arrow{s} \arrow{e}
    \node{\pi_{k+1}} \arrow{e} \arrow{sseee,..}
    \node{\cdots} \arrow[3]{e}
    \node[3]{\varinjlim \pi_k}
    \arrow[2]{s,r}{\hat \phi}
    \\
    \node{G_{n_0}} \arrow{e} \arrow{s,=}
    \node{\cdots} \arrow{e}
    \node{G_{n_k}} \arrow{seee,=}
    \\
    \node{G_0} \arrow{e}
    \node{\cdots} \arrow{e}
    \node{G_k} \arrow{e} \arrow{n}
    \node{G_{k+1}} \arrow{e}
    \node{\cdots} \arrow{e}
    \node{G_{n_k}} \arrow{e}
    \node{\cdots} \arrow{e}
    \node{\varinjlim G_k}
  \end{diagram}
  \]
  Since $\pi_{k+1}$ is finitely generated, the composition $\pi_{k+1}
  \to \varinjlim \pi_k \to \varinjlim G_k$ factors through $G_{r}$ for
  some $r > n_k$.  Note that the two compositions $\pi_k \to \pi_{k+1}
  \to G_r$ and $\pi_k \to G_{n_k} \to G_r$ may not be identical.
  However, composing $G_r \to \varinjlim G_k$ with them, we obtain the
  same maps.  Since $\pi_k$ is finitely generated, it follows that
  $\pi_k \to \pi_{k+1} \to G_{s}$ and $\pi_k \to G_{n_k} \to G_s$ are
  identical for some $s\ge r$.  We choose $s$ as~$n_{k+1}$.  Our
  induction hypothesis is maintained so that the construction of
  $\{n_k\}$ can be continued.

  Now, letting $P_k = G_{n_k}$, obviously $P_k \to P_{k+1}$, $G_k \to
  P_k$ are in $\Omega^R$, and $\varinjlim P_k \cong \varinjlim G_k$.
\end{proof}

\bibliographystyle{amsplainabbrv}
\bibliography{research}

\end{document}